\documentclass[12pt,a4paper,reqno]{amsart}

\usepackage{amsfonts}
\usepackage{amsmath}
\usepackage{amssymb}
\usepackage{amsthm}
\usepackage{amsxtra}
\usepackage{bbm}
\usepackage{braket}
\usepackage{comment}
\usepackage{extarrows}
\usepackage{graphicx}
\usepackage{latexsym}
\usepackage{mathrsfs}
\usepackage{yhmath}
\usepackage{nicematrix}
\usepackage{mathtools}

\usepackage{float}
\usepackage{booktabs}
\usepackage{verbatim}
\usepackage{xcolor}

\usepackage[pdfborder={0 0 0}]{hyperref} 
\hypersetup{breaklinks, raiselinks}

\usepackage{bm}

\usepackage{cleveref}

\usepackage{a4wide}
\usepackage{fullpage}

\theoremstyle{plain}
\newtheorem{Theorem}{Theorem}[section]
\newtheorem{Lemma}[Theorem]{Lemma}
\newtheorem{Corollary}[Theorem]{Corollary}
\newtheorem{Proposition}[Theorem]{Proposition}

\theoremstyle{definition}

\newtheorem{Assumptions and Discussion}[Theorem]{Assumptions and Discussion}

\theoremstyle{remark}

\newtheorem*{acknowledgment*}{Acknowledgment}

\numberwithin{equation}{section}

\usepackage[shortlabels]{enumitem}
\SetEnumerateShortLabel{a}{\textup{(\alph*)}}
\SetEnumerateShortLabel{A}{\textup{(\Alph*)}}
\SetEnumerateShortLabel{1}{\textup{(\arabic*)}}
\SetEnumerateShortLabel{i}{\textup{(\roman*)}}
\SetEnumerateShortLabel{I}{\textup{(\Roman*)}}

\def\bar#1{\overline{#1}}

\begin{document}

\title{On pseudoholomorphic map between almost Hermitian manifolds}

\author[Chiakuei Peng, Xiaowei Xu]{ Chiakuei Peng \quad Xiaowei Xu*}





\address{School of Mathematical Sciences, University of Chinese Academy of Sciences, Beijing, 100049, P.R.~China.}
\email{pengck@ucas.edu.cn}

\address{CAS Wu Wen-Tsun Key Laboratory of Mathematics, School of Mathematical Sciences, University of Science and Technology of China, Hefei, Anhui, 230026, P.R.~China.}
\email{xwxu09@ustc.edu.cn}

\begin{abstract}
In this paper, we use the canonical connection instead of Levi-Civita connection to study the smooth maps between almost Hermitian manifolds, especially, the pseudoholomorphic ones. By using the Bochner formulas, we obtian the $C^2$-estimate of canonical second fundamental form, Liouville type theorems of pseudoholomorphic maps, pseudoholomorphicity of pluriharmonic maps, and Simons integral inequality of pseudoholomorphic isometric immersion.
\end{abstract}
\maketitle
\section{Introduction}
In 1985, M.Gromov (\cite{[G-85]}) introduced the beautiful theory of pseudoholomorphic curves in symplectic manifold, which has profoundly influenced the study of symplectic geometry and symplectic topology. One can refer to D.Macduff and D.Salamon's book \cite{[MS-12]} and references therein. Moreover,  Gromov's compactness theorem of pseudoholomorphic curves was also a interesting probelm from the analytic point of views.
T.H.Parker and J.G.Wolfson (\cite{[PW-93]}), R.G.Ye (\cite{[Y-94]}) independently gave the analytic proofs of Gromov's compactness theorem.
This field has been extended to the case of domain manifold with higher dimention, i.e., the pseudoholomorphic map between almost Hermitian manifolds.
Ch.Y.Wang (\cite{[W1-03]}) studied the regularity and blow-up analysis of pseudoholomorphic maps by using the techniques developed in the theory of harmonic maps. T.Rivi\'ere and G.Tian (\cite{[RT-04]}) studied the almost complex four-manifold into algebraic varieties  in connection with C.Taubes' works. For the triholmorphic maps, Ch,Y.Wang (\cite{[W2-03]}), C.Bellettini and G.Tian (\cite{[BT-19]}) obtian the quantization of energy, compactness results respectively. In this paper, we wish to focus on the geometry of pseudoholomorphic maps by using the canonical connection on almost Hermitian manifold instead of the Levi-Civita connection.

The pseudoholomorphic map is closed related to the theory of harmonic map. It is known that (\cite{[EL-78]}, \cite{[X-96]}) a holomorphic map between K\"ahlerian manifolds is harmonic. Conversely, how to judge a harmonic map to be a holomorphic one is an important topic in the theory of harmonic map.
Y.T.Siu (\cite{[S-80]}) proved the holomorphicity result of harmonic maps between compact K\"ahler manifold under the target manifold has strongly negative curvature. J.Jost and S.T.Yau (\cite{[JY-93]}) used the Hermitian harmonic maps to study maps from Hermitian manifold to Riemannian manifods, and they also obtianed rigidity results in Hermitian geometry. K.F.Liu and X.K. Yang (\cite{[LY-14]}) extended Siu's result to Hermitian manifolds, and they also studied several types of harmonic maps from Hermitian manifold. Y.X.Dong (\cite{[D-13]}) has used the monotonicity formulae to derive holomorphicity and Liouville type results for pluriharmonic maps and harmonic maps between K\"ahlerian manifolds. Y.L.Xin (\cite{[X-96]}) found a sufficient condition of partial energe density to judge the holomorphicity of a harmonic maps from Riemannian surfaces into complex projective space. 
L.Ni (\cite{[N-19]}) use the $\partial\bar{\partial}$-Bochner formulae to study the general Schwarz lemma and their applications on K\"ahler manifolds.
K.Tang (\cite{[T-19]}) generalized works of L.Ni to

\noindent------------------------------------------

*The corresponding author.

\noindent  the case of Hermitian manifolds.
It is clear that the Bochner technique, i.e., the Bochner type formulae of different geometric quantity, palys an important role in studing holomorphicity and Liouville type results of harmonic maps.

The pseudoholomorphic isometry is a special class of pseudoholomorphic map that preserves the metric. If we view from theory of submanifold, the geometry of pseudoholomorphic isometry between almost Hermitain manifolds is very similar to the minimal one between Riemannian manifolds. For examples, one can choose adapted unitary frame field along a pseudoholomorphic isometry, and the trace of canonical second fundamental form is vanishing. In order to solve the Bernstein problem, J.Simons (\cite{[J-68]}) calculated the Laplacian of squre norm of the second fundamental form of closed minimal submanifolds in the unit sphere, and he obtianed the famous Simons integral inequality. K.Ogiue (\cite{[O-74]}) derived the complex version of the formula for holomorphic isometry between K\"ahler manifolds. H.W.Xu (\cite{[X-95]}) generalized Simons integral inequality to closed minimal submanifold in pinched Riemannian manifolds. We wish to develop these works to pseudoholomorphic isometry between almost Hermitian manifolds.

Notice that the canonical connection on almost Hermitian manifold  is the unique one preserves the metric, the almost complex structure and also has vanishing $(1,1)$-part of torsion. So, it's reasonable and natural to use the canonical connetion to study the properties of maps between almost Hermitian manifolds which are independent of Levi-Civita connection, such as  the Liouville type results, the pseudoholomorphicity and so on. Inspired by the works of V.Tosatti (\cite{[T-07]}), V.Tosatti, B.Wenkove and S.T.Yau (\cite{[TWY-08]}), X.Zhang (\cite{[Zh-12]}), by using the canonical connection, one can find the Bochner formulae in complex version is quite like the one of harmonic map. We utilize techniques and methods of harmonic maps in the nice book (\cite{[X-96]}) by Y.L.Xin to study the pseudoholomorphic maps  and  pluriharmonic maps. These are the main motivation of our paper.

The paper organized as follows. In section 2, we recall the almost Hermitian geometry with respect to the canonical connection. In section 3, we calculate the canonical Laplacian of partial energy density of pluriharmonic maps, and Laplacian of the half norm of canonical second fundamental form $S$ of pseudoholomorphic maps. In section 4 and 5, we study the Liouville type results of pseudoholomorphic  maps and pseudoholomorphicity of pluriharminic maps by using the Bochner type formulae obtained, respectively. In section 6, we derivie Simons integral inequality for pseudoholomorphic isometry, and we also give the bound of $S$ under the condition of parallel canonical second fundamental form.

\section{Geometries of almost Hermitian manifold}
The purpose of this section is to establish notation, review some well known facts (see \cite{[T-07]}, \cite{[TWY-08]}) concerning almost Hermitian geometry, and to establish the basic identities of smooth map between almost Hermitian manifolds related to the canonical connection. For the curvature properties, one can refer to the detailed papers \cite{[FZ-19]}, \cite{[Y-17]} and references therein.

\subsection{Geometry of canonical connection}

Let $(M, J, g)$ be a $2m$-dimensional almost Hermitian manifold, which means that $J$ is an almost complex structure on $M$ and Riemannian metric $g$ is $J$-invariant. 
The tangent bundle of $M$ denoted by $TM$, and its complexification $TM\otimes\mathbb{C}$ denoted by $T^\mathbb{C}M$. By extending the almost complex structure 
$\mathbb{C}$-linearly to $T^\mathbb{C}M$, we obtain the decomposition
\begin{equation}\label{2-1-1}
T^\mathbb{C}M=T'M\oplus T''M,
\end{equation}
where $T'M$ and $T''M$ are the eigenspaces of $J$ corresponding to eigenvalues $\sqrt{-1}$ and $-\sqrt{-1}$ respectively. Similarly, by extending $J$ to forms, one can decompose $k$-forms into a sum of $(p,q)$-forms with $p+q=k$. Extending the almost Hermitian metric $\mathbb{C}$-linearly to $T^\mathbb{C}M$, then it induces a Hermitian metric on $T'M$, still denoted by $g$. Locally, we can choose a unitary frame $\{e_1,\ldots,e_m\}$ with the dual coframe $\{\theta^1,\ldots,\theta^m\}$ so that
\begin{equation}\label{2-1-2}
g=\theta^i\otimes\overline{\theta^i},\hspace{1cm} dV_g=(\sqrt{-1})^m\theta^1\wedge\overline{\theta^1}\wedge\cdots\wedge\theta^m\wedge\overline{\theta^m}.
\end{equation}
The fundamental 2-form corresponding to the Hermitian metric $g$ is given by
\begin{equation}\label{2-1-3}
\omega=\sqrt{-1}\;\theta^i\wedge\overline{\theta^i}.
\end{equation}

Let $\nabla$ be an affine connection on $TM$, which is called \emph{almost-Hermitian} if it satifies
\begin{equation}\label{2-1-4}
\nabla J=0,\hspace{1cm} \nabla g=0.
\end{equation}
It is easy to see that such connections always exist on any almost Hermitian manifold. One can refer to \cite{[FZ-19]}, \cite{[G-97]} for details.  Locally, it follows from (\ref{2-1-4}) that the \emph{connection $1$-forms} $\{\theta^j_i\}$ are given by
\begin{equation}\label{2-1-5}
\nabla e_i=\theta_i^j e_j,\hspace{1cm}\theta^j_i+\overline{\theta^i_j}=0.
\end{equation}
The first and second structure equations of $\nabla$ are determined by
\begin{equation}\label{2-1-6}
d\theta^i=-\theta^i_j\wedge\theta^j+\Theta^i,
\end{equation}
\begin{equation}\label{2-1-7}
d\theta^i_j=-\theta^i_k\wedge\theta^k_j+\Omega^i_j,\hspace{0.5cm}\Omega^j_i+\overline{\Omega^i_j}=0,
\end{equation}
where the $2$-forms $\Theta^i$, $\Omega^i_j$ are called \emph{torsion} and \emph{curvature} of $\nabla$, respectively.
It is well known (see \cite{[G-97]}) that there exists a unique almost Hermitian connection $\nabla$ on $TM$ whose $(1,1)$-part satisfies
\begin{equation}\label{2-1-8}
(\Theta^i)^{(1,1)}=0.
\end{equation}
which is called \emph{canonical connection} on $(M,J,g)$. It is just the Chern connection when $J$ is integrable. We will always to use $\nabla$ to stand for the canonical connection on almost Hermitian manifold.


There are several important families of almost Hermitian manifold defined by using the fundamental 2-form, Levi-Civita connection and almost complex structure. We say that $(M,J,g)$ is \emph{almost-K\"ahlerian} if $d\omega=0$, that it is \emph{quasi-K\"ahlerian} if $(d\omega)^{(1,2)}=0$, and that it is \emph{semi-K\"ahlerian} if $d^*\omega=0$, where $d^*$ is the codifferential operator of $d$ with respect to the metric $g$. To character them, we set
\begin{equation}\label{2-1-11}
\Theta^i:=L^i_{jk}\,\theta^j\wedge\theta^k+N^i_{\bar{j}\,\bar{k}}\,\overline{\theta^j}\wedge\overline{\theta^k}
\end{equation}
with $L^i_{jk}=-L^i_{kj}$, $N^i_{\bar{j}\,\bar{k}}=-N^i_{\bar{k}\,\bar{j}}$. Then, in terms of $L^i_{jk}$ and $N^i_{\bar{j}\,\bar{k}}$, it is easy to check that $(M,J,g)$ is almost-K\"ahlerian if and only if
\begin{equation}\label{2-1-12}
L^i_{jk}=0\hspace{0.3cm} and \hspace{0.3cm} N^i_{\bar{j}\,\bar{k}}+N^j_{\bar{k}\,\bar{i}}+N^k_{\bar{i}\,\bar{j}}=0,
\end{equation}
it is quasi-K\"ahlerian if and only if
\begin{equation}\label{2-1-13}
L^i_{jk}=0,
\end{equation}
and it is semi-K\"ahlerian if and only if
\begin{equation}\label{2-1-14}
L_j:=\sum\limits_{i}L^i_{ij}=0.
\end{equation}
On the other hand, by using the Levi-Civita connection, we call $(M,J,g)$ is \emph{nearly-K\"ahlerian} if $(D_XJ)X=0$ for any tangent vector $X$ on $M$. This is equivalent to
\begin{equation}\label{2-2-60}
L^i_{jk}=0\hspace{0.3cm}and \hspace{0.3cm} N^i_{\bar{j}\,\bar{k}}=N^j_{\bar{k}\,\bar{i}}.
\end{equation}

We also set
\begin{equation}\label{2-1-15}
\Omega^i_j:=R^i_{jk\ell}\,\theta^k\wedge\theta^\ell+R^i_{jk\bar{\ell}}\,\theta^k\wedge\overline{\theta^\ell}+R^i_{j\bar{k}\,\bar{\ell}}\,\overline{\theta^k}\wedge\overline{\theta^\ell}
\end{equation}
with $R^i_{jk\ell}=-R^i_{j\ell k}$, $R^i_{j\bar{k}\,\bar{\ell}}=-R^i_{j\bar{\ell}\,\bar{k}}$. The skew-Hermitian of curvature forms $\Omega^i_j$ implies
\begin{equation}\label{2-1-16}
R^i_{jk\ell}=\overline{R^j_{i\bar{\ell}\,\bar{k}}},\hspace{1cm}R^i_{jk\bar{\ell}}=\overline{R^j_{i\ell\bar{k}}}.
\end{equation}
For a nonzero vectors $X\in T'M$, we call $HS(X):=R^j_{ik\bar{\ell}}\,X^i\overline{X^j}X^k\overline{X^\ell}/|X|^4$ the \emph{holomorphic sectional curvature} of canonical connection in direction $X$, and
for two nonzero vectors $X,Y\in T'M$ we call
$HB(X,Y):=R^j_{ik\bar{\ell}}\,X^i\overline{X^j}Y^k\overline{Y^\ell}/(|X|^2|Y|^2)$ the \emph{holomophic bisectional curvature} of canonical connection in the directions $X$, $Y$. We say the holomorphic bisectional curvature is bounded below by $A$ and bounded above by $B$, if $A\leq HB(X,Y)\leq B$ for all nonzero $(1,0)$-type vectors $X,Y$. The tensors $R'_{k\bar{\ell}}:=R^i_{ik\bar{\ell}}$, $R''_{k\bar{\ell}}:=R^\ell_{ki\bar{i}}$ are called the \emph{first Ricci curvature}, \emph{second Ricci curvature} of the canonical connection, which will be denoted by $Ric'_M$, $Ric''_M$, respectively. We say the first (resp. second) Ricci curvature is bounded below by $A$ if $R'_{k\bar{\ell}}X^k\overline{X^\ell}\geq A|X|^2$ (resp. $R''_{k\bar{\ell}}X^k\overline{X^\ell}\geq A|X|^2$) for all $(1,0)$-type vector $X$. Similar bounds of torsion and $(2,0)$-part of curvature can be defined. We will say $(M,J,g)$ has \emph{bouded geometry} if
its curvature, convariant derivatives of curvature are uniformly bounded.

\vspace{0.2cm}

Let $u$ be a smooth function on $M$, its differential can be written as
\begin{equation}\label{2-1-32}
du=u_i\,\theta^i+u_{\bar{i}}\,\overline{\theta^i},
\end{equation}
where $u_{\bar{i}}=\overline{u_i}$. Taking the exterior derivative of (\ref{2-1-32}) and using (\ref{2-1-6}), we obtain
\begin{equation}\label{2-1-33}
(u_{i,j}\,\theta^j+u_{i,\bar{j}}\,\overline{\theta^j})\wedge\theta^i+u_i\,\Theta^i+(u_{\bar{i},j}\,\theta^j+u_{\bar{i},\bar{j}}\,\overline{\theta^j})\wedge\overline{\theta^i}+u_{\bar{i}}\,\overline{\Theta^i}=0,
\end{equation}
where $u_{i,j}$, $u_{i,\bar{j}}$, $u_{\bar{i},j}$, $u_{\bar{i},\bar{j}}$ are defined by
\begin{equation}\label{2-1-34}
u_{i,j}\,\theta^j+u_{i,\bar{j}}\,\overline{\theta^j}:=du_i-u_j\,\theta^j_i,
\end{equation}
\begin{equation}\label{2-1-35}
u_{\bar{i},j}\,\theta^j+u_{\bar{i},\bar{j}}\,\overline{\theta^j}:=du_{\bar{i}}-u_{\bar{j}}\,\overline{\theta^j_i}.
\end{equation}
Comparing the coefficients of (\ref{2-1-33}), we obtain
\begin{equation}\label{2-1-36}
u_{i,j}-u_{j,i}=2u_kL^k_{ij}+2u_{\bar{k}}\overline{N^k_{\bar{i}\,\bar{j}}},
\end{equation}
\begin{equation}\label{2-1-37}
u_{i,\bar{j}}-u_{\bar{j},i}=0.
\end{equation}
The \emph{canonical Hessian} of a smooth function $u$ on $M$ is defined by $\nabla du$, which can be expressed as
\begin{equation}\label{2-1-44}
\nabla du=u_{i,j}\,\theta^j\otimes\theta^i+u_{i,\bar{j}}\,\overline{\theta^j}\otimes\theta^i+u_{\bar{i},j}\,\theta^j\otimes\overline{\theta^i}
+u_{\bar{i},\bar{j}}\,\overline{\theta^j}\otimes\overline{\theta^i}.
\end{equation}
By taking the trace of (\ref{2-1-44}), the \emph{canonical Laplacian} of $u$ is defined by
\begin{equation}\label{2-1-45}
\Box u:=\sum\limits_{i}(u_{i,\bar{i}}+u_{\bar{i},i})=2\sum\limits_{i}u_{i,\bar{i}}.
\end{equation}
We call a smooth function $u$ is \emph{strictly plurisubharmonic} if $\nabla du(X,\overline{X})=u_{i,\bar{j}}X^i\overline{X^j}>0$ for any nonzero vector $X\in T'M$.

Let $\Delta$ be the usual Laplacian related to Levi-Civita connection of $g$. The difference (see \cite{[T-07]}, \cite{[TWY-08]}) between $\Delta$ and $\Box$ are given by
\begin{equation}\label{2-2-62}
\Delta u=\Box u-2\langle\nabla u,X_L\rangle,
\end{equation}
for a function $u\in C^2(M)$,
where $\nabla u=u_{\bar{i}}\,e_i+u_i\,\overline{e_i}$, $X_L=\overline{L_i}\,e_i+L_i\,\overline{e_i}$ and $L_i$ defined in (\ref{2-1-14}).
Similarly, for a vector field $X=X^k\,e_k+\overline{X^k}\,\overline{e_k}$, by using the canonical connection and Levi-Civita connection, one can define the \emph{canonical divergence} and usual divergence of $X$, denoted by $\mbox{div}^c(X)$ and $\mbox{div}(X)$ respectively. It is easy to check that the difference of such two divergences is given by
\begin{equation}\label{2-2-62-1}
\mbox{div}(X)=\mbox{div}^c(X)-2\langle X,X_L\rangle.
\end{equation}
Thus, these two Laplacian and two divergences are equal respectively, when $(M,J,g)$ is semi-K\"ahlerian, or quasi-K\"ahlerian, or almost-K\"ahlerian, or nearly-K\"ahlerian.

\subsection{Smooth map between almost Hermitian manifolds} 
Let $(M,J,g)$ and $(\widetilde{M},\tilde{J}, \tilde{g})$ be two almost Hermitiant manifolds of dimension $2m$ and $2n$, with the canonical connections $\nabla$, $\widetilde{\nabla}$, respectively.
Let $f$ be a smooth map from $M$ into $\widetilde{M}$.
Locally, choosing a unitary frame $\{e_i\}$ with the dual coframe $\{\theta^i\}$ on $M$, and choosing a unitary frame $\{\tilde{e}_\alpha\}$ with the dual coframe $\{\tilde{\theta}^\alpha\}$ on $\widetilde{M}$. Set
\begin{equation}\label{2-3-69}
f^*\tilde{\theta}^\alpha:=a^\alpha_i\,\theta^i+a^\alpha_{\bar{i}}\,\overline{\theta^i},\hspace{0.5cm}
f^*\overline{\tilde{\theta}^\alpha}:=a^{\bar{\alpha}}_i\,\theta^i+a^{\bar{\alpha}}_{\bar{i}}\,\overline{\theta^i},
\end{equation}
where $a^\alpha_i$, $a^\alpha_{\bar{i}}$, $a^{\bar{\alpha}}_i$, $a^{\bar{\alpha}}_{\bar{i}}$ satisfy
\begin{equation}\label{2-3-70}
\overline{a^\alpha_i}=a^{\bar{\alpha}}_{\bar{i}}, \hspace{1cm}\overline{a^\alpha_{\bar{i}}}=a^{\bar{\alpha}}_i.
\end{equation}

Taking the exterior derivative of the first identity in (\ref{2-3-69}), using the first structure equations of $\nabla$ and $\widetilde{\nabla}$, we obtain
\begin{equation}\label{2-3-72}
(a^\alpha_{i,j}\,\theta^j+a^\alpha_{i,\bar{j}}\,\overline{\theta}^j)\wedge\theta^i+(a^\alpha_{\bar{i},j}\,\theta^j+a^\alpha_{\bar{i},\bar{j}}\,\overline{\theta}^j)\wedge\overline{\theta^i}+a^\alpha_p\,\Theta^p+a^\alpha_{\bar{p}}\,\overline{\Theta^p}-\widetilde{\Theta}^\alpha=0,
\end{equation}
where $a^\alpha_{i,j}$, $a^\alpha_{i,\bar{j}}$, $a^\alpha_{\bar{i},j}$,  $a^\alpha_{\bar{i},\bar{j}}$ are defined by
\begin{equation}\label{2-3-73}
a^\alpha_{i,j}\,\theta^j+a^\alpha_{i,\bar{j}}\,\overline{\theta^j}:=da^\alpha_i-a^\alpha_j\,\theta^j_i+a^\beta_i\,\tilde{\theta}^\alpha_\beta,
\end{equation}
\begin{equation}\label{2-3-74}
a^\alpha_{\bar{i},j}\,\theta^j+a^\alpha_{\bar{i},\bar{j}}\,\overline{\theta^j}:=da^\alpha_{\bar{i}}-a^\alpha_{\bar{j}}\,\overline{\theta^j_i}+a^\beta_{\bar{i}}\,\tilde{\theta}^\alpha_\beta.
\end{equation}
Here we have omitted the pull-back $f^*$ acts on  $\widetilde{\Theta}^\alpha$, $\tilde{\theta}^\alpha_\beta$, and similar conventions will be used in the sequal. Comparing the coefficients of (\ref{2-3-72}), we obtain
\begin{equation}\label{2-3-75}
a^\alpha_{i,j}-a^\alpha_{j,i}=2a^\alpha_p\,L^p_{ij}+2a^\alpha_{\bar{p}}\,\overline{N^p_{\bar{i}\,\bar{j}}}-2a^\beta_ia^\gamma_j\tilde{L}^\alpha_{\beta\gamma}-2\overline{a^\beta_{\bar{i}}}\,\overline{a^\gamma_{\bar{j}}}\,\widetilde{N}^\alpha_{\bar{\beta}\bar{\gamma}},
\end{equation}
\begin{equation}\label{2-3-76}
a^\alpha_{i,\bar{j}}-a^\alpha_{\bar{j},i}=(a^\beta_{\bar{j}}a^\gamma_i-a^\beta_ia^\gamma_{\bar{j}})\tilde{L}^\alpha_{\beta\gamma}+(\overline{a^\beta_j}\,\overline{a^\gamma_{\bar{i}}}-\overline{a^\beta_{\bar{i}}}\,\overline{a^\gamma_j})\widetilde{N}^\alpha_{\bar{\beta}\bar{\gamma}},
\end{equation}
\begin{equation}\label{2-3-77}
a^\alpha_{\bar{i},\bar{j}}-a^\alpha_{\bar{j},\bar{i}}=2a^\alpha_p\,N^p_{\bar{i}\,\bar{j}}+2a^\alpha_{\bar{p}}\,\overline{L^p_{ij}}-2a^\beta_{\bar{i}}a^\gamma_{\bar{j}}\tilde{L}^\alpha_{\beta\gamma}-2\overline{a^\beta_{i}}\,\overline{a^\gamma_{j}}\,\widetilde{N}^\alpha_{\bar{\beta}\bar{\gamma}}.
\end{equation}

Taking the exterior differential of (\ref{2-3-73}), one have
\begin{equation}\label{2-3-78}
(a^\alpha_{i,jk}\,\theta^k+a^\alpha_{i,j\bar{k}}\,\overline{\theta^k})\wedge\theta^j+(a^\alpha_{i,\bar{j}k}\,\theta^k+a^\alpha_{i,\bar{j}\bar{k}}\,\overline{\theta^k})\wedge\overline{\theta^j}+a^\alpha_{i,p}\Theta^p+a^\alpha_{i,\bar{p}}\overline{\Theta^p}=-a^\alpha_p\,\Omega^p_i+a^\beta_i\widetilde{\Omega}^\alpha_\beta,
\end{equation}
which implies
\begin{eqnarray}\label{2-3-79}
a^\alpha_{i,jk}-a^\alpha_{i,kj}&=&2a^\alpha_p\,R^p_{ijk}-2a^\beta_ia^\gamma_ja^\delta_k\widetilde{R}^\alpha_{\beta\gamma\delta}-a^\beta_i(a^\gamma_ja^{\bar{\delta}}_k-a^\gamma_ka^{\bar{\delta}}_j)\widetilde{R}^\alpha_{\beta\gamma\bar{\delta}}\nonumber\\
&{}&-2a^\beta_ia^{\bar{\gamma}}_ja^{\bar{\delta}}_k\widetilde{R}^\alpha_{\beta\bar{\gamma}\bar{\delta}}+2a^\alpha_{i,p}L^p_{jk}+2a^\alpha_{i,\bar{p}}\,\overline{N^p_{\bar{j}\,\bar{k}}},
\end{eqnarray}
\begin{eqnarray}\label{2-3-80}
a^\alpha_{i,j\bar{k}}-a^\alpha_{i,\bar{k}j}=a^\alpha_pR^p_{ij\bar{k}}-2a^\beta_ia^\gamma_ja^\delta_{\bar{k}}\widetilde{R}^\alpha_{\beta\gamma\delta}
-a^\beta_i(a^\gamma_ja^{\bar{\delta}}_{\bar{k}}-a^\gamma_{\bar{k}}a^{\bar{\delta}}_{j})\widetilde{R}^\alpha_{\beta\gamma\bar{\delta}}-2a^\beta_ia^{\bar{\gamma}}_ja^{\bar{\delta}}_{\bar{k}}\widetilde{R}^\alpha_{\beta\bar{\gamma}\bar{\delta}},
\end{eqnarray}
\begin{eqnarray}\label{2-3-81}
a^\alpha_{i,\bar{j}\,\bar{k}}-a^\alpha_{i,\bar{k}\,\bar{j}}&=&2a^\alpha_p\,R^p_{i\bar{j}\,\bar{k}}-2a^\beta_ia^\gamma_{\bar{j}}a^\delta_{\bar{k}}\widetilde{R}^\alpha_{\beta\gamma\delta}-a^\beta_i(a^\gamma_{\bar{j}}a^{\bar{\delta}}_{\bar{k}}-a^\gamma_{\bar{k}}a^{\bar{\delta}}_{\bar{j}})\widetilde{R}^\alpha_{\beta\gamma\bar{\delta}}\nonumber\\
&{}&-2a^\beta_ia^{\bar{\gamma}}_{\bar{j}}a^{\bar{\delta}}_{\bar{k}}\widetilde{R}^\alpha_{\beta\bar{\gamma}\bar{\delta}}+2a^\alpha_{i,p}N^p_{\bar{j}\,\bar{k}}+2a^\alpha_{i,\bar{p}}\,\overline{L^p_{jk}},
\end{eqnarray}
where $a^\alpha_{i,jk}$, $a^\alpha_{i,j\bar{k}}$, $a^\alpha_{i,\bar{j}k}$, $a^\alpha_{i,\bar{j}\bar{k}}$ defined by
\begin{equation}\label{2-3-82}
a^\alpha_{i,jk}\,\theta^k+a^\alpha_{i,j\bar{k}}\,\overline{\theta^k}:=da^\alpha_{i,j}-a^\alpha_{p,j}\,\theta^p_i-a^\alpha_{i,p}\,\theta^p_j+a^\beta_{i,j}\tilde{\theta}_\beta^\alpha,
\end{equation}
\begin{equation}\label{2-3-83}
a^\alpha_{i,\bar{j}k}\,\theta^k+a^\alpha_{i,\bar{j}\,\bar{k}}\,\overline{\theta^k}:=da^\alpha_{i,\bar{j}}-a^\alpha_{p,\bar{j}}\,\theta^p_i-a^\alpha_{i,\bar{p}}\,\overline{\theta^p_j}+a^\beta_{i,\bar{j}}\tilde{\theta}_\beta^\alpha.
\end{equation}

Taking the exterior derivative of (\ref{2-3-74}), we have
\begin{equation}\label{2-3-83-1}
(a^\alpha_{\bar{i},jk}\theta^k+a^\alpha_{\bar{i},j\bar{k}}\overline{\theta^k})\wedge\theta^j
+(a^\alpha_{\bar{i},\bar{j}k}\theta^k+a^\alpha_{\bar{i},\bar{j}\,\bar{k}}\overline{\theta^k})\wedge\overline{\theta^j}
+a^\alpha_{\bar{i},p}\Theta^p+a^\alpha_{\bar{i},\bar{p}}\overline{\Theta^p}
=-a^\alpha_{\bar{p}}\overline{\Omega^p_i}+a^\beta_{\bar{i}}\widetilde{\Omega}^\alpha_\beta,
\end{equation}
which implies
\begin{eqnarray}\label{2-3-83-2}
a^\alpha_{\bar{i},jk}-a^\alpha_{\bar{i},kj}&=&2a^\alpha_{\bar{p}}\overline{R^p_{i\bar{j}\bar{k}}}-2a^\beta_{\bar{i}}a^\gamma_ja^\delta_k\widetilde{R}^\alpha_{\beta\gamma\delta}-a^\beta_{\bar{i}}(a^\gamma_ja^{\bar{\delta}}_k-a^\gamma_ka^{\bar{\delta}}_j)\widetilde{R}^\alpha_{\beta\gamma\bar{\delta}}\nonumber\\
&{}&-2a^\beta_{\bar{i}}a^{\bar{\gamma}}_ja^{\bar{\delta}}_k\widetilde{R}^\alpha_{\beta\bar{\gamma}\bar{\delta}}+2a^\alpha_{\bar{i}p}L^p_{jk}+2a^\alpha_{\bar{i}\bar{p}}\overline{N^p_{\bar{j}\bar{k}}},
\end{eqnarray}
\begin{eqnarray}\label{2-3-83-3}
a^\alpha_{\bar{i},j\bar{k}}-a^\alpha_{\bar{i},\bar{k}j}=-a^\alpha_{\bar{p}}\overline{R^p_{ik\bar{j}}}
-2a^\beta_{\bar{i}}a^\gamma_ja^\delta_{\bar{k}}\widetilde{R}^\alpha_{\beta\gamma\delta}
-a^\beta_{\bar{i}}(a^\gamma_ja^{\bar{\delta}}_{\bar{k}}-a^\gamma_{\bar{k}}a^{\bar{\delta}}_j)\widetilde{R}^\alpha_{\beta\gamma\bar{\delta}}
-2a^\beta_{\bar{i}}a^{\bar{\gamma}}_ja^{\bar{\delta}}_{\bar{k}}\widetilde{R}^\alpha_{\beta\bar{\gamma}\bar{\delta}},
\end{eqnarray}
\begin{eqnarray}\label{2-3-83-4}
a^\alpha_{\bar{i},\bar{j}\,\bar{k}}-a^\alpha_{\bar{i},\bar{k}\,\bar{j}}&=&2a^\alpha_{\bar{p}}\overline{R^p_{ijk}}-2a^\beta_{\bar{i}}a^\gamma_{\bar{j}}a^\delta_{\bar{k}}\widetilde{R}^\alpha_{\beta\gamma\delta}-a^\beta_{\bar{i}}(a^\gamma_{\bar{j}}a^{\bar{\delta}}_{\bar{k}}-a^\gamma_{\bar{k}}a^{\bar{\delta}}_{\bar{j}})\widetilde{R}^\alpha_{\beta\gamma\bar{\delta}}\nonumber\\
&{}&-2a^\beta_{\bar{i}}a^{\bar{\gamma}}_{\bar{j}}a^{\bar{\delta}}_{\bar{k}}\widetilde{R}^\alpha_{\beta\bar{\gamma}\bar{\delta}}+2a^\alpha_{\bar{i}p}N^p_{\bar{j}\bar{k}}+2a^\alpha_{\bar{i}\bar{p}}\overline{L^p_{jk}},
\end{eqnarray}
where $a^\alpha_{\bar{i},jk}$, $a^\alpha_{\bar{i},j\bar{k}}$, $a^\alpha_{\bar{i},\bar{j}k}$ and $a^\alpha_{\bar{i},\bar{j}\bar{k}}$ defined by
\begin{equation}\label{2-3-83-5}
a^\alpha_{\bar{i},jk}\theta^k+a^\alpha_{\bar{i},j\bar{k}}\,\overline{\theta^k}:=da^\alpha_{\bar{i},j}-a^\alpha_{\bar{k},j}\overline{\theta^k_i}
-a^\alpha_{\bar{i},k}\theta^k_j+a^\beta_{\bar{i},j}\tilde{\theta}_\beta^\alpha,
\end{equation}
\begin{equation}\label{2-3-83-6}
a^\alpha_{\bar{i},\bar{j}k}\theta^k+a^\alpha_{\bar{i},\bar{j}\,\bar{k}}\,\overline{\theta^k}:=da^\alpha_{\bar{i},\bar{j}}-a^\alpha_{\bar{k},\bar{j}}\overline{\theta^k_i}
-a^\alpha_{\bar{i},\bar{k}}\overline{\theta^k_j}+a^\beta_{\bar{i},\bar{j}}\tilde{\theta}_\beta^\alpha.
\end{equation}
Similarly, taking the exterior derivative of (\ref{2-3-82}), by comparing the (1,1)-part, one can get
\begin{eqnarray}\label{2-3-84}
a^\alpha_{i,jk\bar{\ell}}-a^\alpha_{i,j\bar{\ell}k}&=&a^\alpha_{p,j}R^p_{ik\bar{\ell}}+a^\alpha_{i,p}R^p_{jk\bar{\ell}}-2a^\beta_{i,j}a^\gamma_ka^\delta_{\bar{\ell}}\widetilde{R}^\alpha_{\beta\gamma\delta}-a^\beta_{i,j}(a^\gamma_ka^{\bar{\delta}}_{\bar{\ell}}-a^\gamma_{\bar{\ell}}a^{\bar{\delta}}_k)\widetilde{R}^\alpha_{\beta\gamma\bar{\delta}}\nonumber\\
&{}&-2a^\beta_{i,j}a^{\bar{\gamma}}_ka^{\bar{\delta}}_{\bar{\ell}}\widetilde{R}^\alpha_{\beta\bar{\gamma}\,\bar{\delta}},
\end{eqnarray}
where the covariant derivatives $a^\alpha_{i,jk\bar{\ell}}$, $a^\alpha_{i,j\bar{\ell}k}$ defined in the natural way as before.

In terms of unitary frame, from (\ref{2-3-69}), the differential of $f$ can be expressed by
\begin{equation}
df=a^\alpha_i\,\theta^i\otimes\tilde{e}_\alpha+a^\alpha_{\bar{i}}\,\overline{\theta^i}\otimes\tilde{e}_\alpha
+a^{\bar{\alpha}}_i\,\theta^i\otimes\overline{\tilde{e}_\alpha}+a^{\bar{\alpha}}_{\bar{i}}\,\overline{\theta^i}\otimes\overline{\tilde{e}_\alpha},
\end{equation}
which is a smooth section of bundle $T^*M\otimes f^{-1}T\widetilde{M}$ or $(T^{\mathbb{C}}M)^*\otimes f^{-1}T^{\mathbb{C}}\widetilde{M}$. Taking the covariant differntial of $df$ by using the connection induced by the canonical connection $\nabla$ and $\widetilde{\nabla}$, one can get the \emph{canonical second fundamental form} $\nabla df$ of $f$. More explicitly, we have
\begin{eqnarray}
\nabla df&=&a^\alpha_{i,j}\,\theta^j\otimes\theta^i\otimes\tilde{e}_\alpha+a^\alpha_{i,\bar{j}}\,\overline{\theta^j}\otimes\theta^i\otimes\tilde{e}_\alpha
+a^\alpha_{\bar{i},j}\,\theta^j\otimes\overline{\theta^i}\otimes\tilde{e}_\alpha+a^\alpha_{\bar{i},\bar{j}}\,\overline{\theta^j}\otimes\overline{\theta^i}\otimes\tilde{e}_\alpha\hspace{0.8cm}\nonumber\\
&+&a^{\bar{\alpha}}_{i,j}\,\theta^j\otimes\theta^i\otimes\overline{\tilde{e}_\alpha}\!+a^{\bar{\alpha}}_{i,\bar{j}}\,\overline{\theta^j}\otimes\theta^i\otimes\overline{\tilde{e}_\alpha}
\!+a^{\bar{\alpha}}_{\bar{i},j}\,\theta^j\otimes\overline{\theta^i}\otimes\overline{\tilde{e}_\alpha}\!+a^{\bar{\alpha}}_{\bar{i},\bar{j}}\,\overline{\theta^j}\otimes\overline{\theta^i}\otimes\overline{\tilde{e}_\alpha}.
\end{eqnarray}
Such a map is called \emph{canonical geodesic} if $\nabla df=0$; and it is called \emph{pluriharmonic} if $(\nabla df)^{(1,1)}=0$, which means that $a^\alpha_{i,\bar{j}}=a^\alpha_{\bar{i},j}=a^{\bar{\alpha}}_{i,\bar{j}}=a^{\bar{\alpha}}_{\bar{i},j}=0$.

The \emph{partial energy densities} $e'(f)$, $e''(f)$, \emph{energy density} $e(f)$ (see \cite{[EL-78]}, \cite{[X-96]}) of a smooth map $f$ between almost Hermitian manifolds defined by
\begin{equation}\label{2-3-86}
e'(f):=\sum\limits_{i,\alpha}|a^\alpha_i|^2,\hspace{0.5cm}
e''(f):=\sum\limits_{i,\alpha}|a^\alpha_{\bar{i}}|^2,\hspace{0.5cm}
e(f):=e'(f)+e''(f),
\end{equation}
and the corresponding \emph{partial energies} and \emph{energy} of $f$ defined by
\begin{equation}\label{2-3-87}
E'(f):=\int_{M}e'(f)\,dV_g,\hspace{0.5cm}
E''(f):=\int_{M}e''(f)\,dV_g,\hspace{0.5cm}
E(f):=\int_{M}e(f)\,dV_g,
\end{equation}
respectively.
A smooth map is called \emph{harmonic} if it is a critical point of the energy functional $E$. For more details about harmonic map one can refer to the books \cite{[EL-78]}, \cite{[X-96]} and references therein.

Let $f$ be a smooth map between $(M,J,g)$ and $(\widetilde{M},\tilde{J}, \tilde{g})$, it is called \emph{pseudoholomorphic} if its differential $df$ satisfies
\begin{equation}\label{2.49}
df J=\tilde{J}df.
\end{equation}
 In terms of differential forms, one can check that the indentity (\ref{2.49}) is equivalent to 
\begin{equation}\label{2-3-71}
a^\alpha_{\bar{i}}=a^{\bar{\alpha}}_i=0,
\end{equation}
which means that $f^*$ preserves the forms of $(1,0)$-type. Thus, the definition (\ref{2-3-74}) and the identity (\ref{2-3-76}) imply
\begin{equation}\label{2.52}
a^\alpha_{\bar{i},j}=a^\alpha_{\bar{i},\bar{j}}=a^\alpha_{i,\bar{j}}=0.
\end{equation}
This tell us that a pseudoholomorphic map is naturally a pluriharmonic one.

\section{Bochner type formulae}

In this section we wish to calculate the canonical Laplacian of the partial energy densities of pluriharmonic maps, and the canonical Laplacian of half norm of canonical second fundamental form of pseudoholomorphic map. 
\begin{Proposition}
Let $f$ be a pluriharmonic map between almost Hermitian manifolds $(M,J,g)$ and $(\widetilde{M},\tilde{J}, \tilde{g})$, then we have

\begin{eqnarray*}
(\,1\,) \;\frac{1}{2}\Box \,e'=|a^\alpha_{i,j}|^2+\overline{a^\alpha_i}a^\alpha_jR''_{i\bar{j}}-4\emph{Re}(\overline{a^\alpha_i}a^\beta_ia^\gamma_ja^\delta_{\bar{j}}\widetilde{R}^\alpha_{\beta\gamma\delta})-\overline{a^\alpha_i}a^\beta_i(a^\gamma_j\overline{a^\delta_j}-a^\gamma_{\bar{j}}\overline{a^\delta_{\bar{j}}})\widetilde{R}^\alpha_{\beta\gamma\bar{\delta}};
\end{eqnarray*}
\begin{eqnarray*}
(\,2\,) \;\frac{1}{2}\Box \,e''=|a^\alpha_{\bar{i},\bar{j}}|^2+a^\alpha_{\bar{i}}\overline{a^\alpha_{\bar{j}}}R''_{i\bar{j}}+4\emph{Re}(a^\alpha_{\bar{i}}\overline{a^\beta_{\bar{i}}}a^\gamma_{\bar{j}}a^\delta_{j}\widetilde{R}^\alpha_{\beta\gamma\delta})+a^\alpha_{\bar{i}}\overline{a^\beta_{\bar{i}}}(\overline{a^\gamma_j}a^\delta_j-\overline{a^\gamma_{\bar{j}}}a^\delta_{\bar{j}})\widetilde{R}^\alpha_{\beta\gamma\bar{\delta}}.
\end{eqnarray*}
\end{Proposition}

\emph{Proof}. Since $f$ is pluriharmonic, i.e., $a^\alpha_{i,\bar{j}}=a^\alpha_{\bar{i},j}=0$, we have
\begin{eqnarray}\label{3-1}
\frac{1}{2}\Box\,e'&=&(a^\alpha_i\,\overline{a^\alpha_i}),_{j\bar{j}}\nonumber\\
&=&(a^\alpha_{i,j}\,\overline{a^\alpha_i}+a^\alpha_i\,\overline{a^\alpha_{i,\bar{j}}}),_{\bar{j}}\nonumber\\
&=&a^\alpha_{i,j\bar{j}}\,\overline{a^\alpha_i}+a^\alpha_{i,j}\,\overline{a^\alpha_{i,j}}.
\end{eqnarray}
By using (\ref{2-3-80}) with $a^\alpha_{i,\bar{j}j}=0$, we obtain
\begin{eqnarray}\label{3-2}
a^\alpha_{i,j\bar{j}}&=&a^\alpha_pR^p_{ij\bar{j}}-2a^\beta_ia^\gamma_ja^\delta_{\bar{j}}\widetilde{R}^\alpha_{\beta\gamma\delta}-a^\beta_i(a^\gamma_j\overline{a^\delta_j}-a^\gamma_{\bar{j}}\overline{a^\delta_{\bar{j}}})\widetilde{R}^\alpha_{\beta\gamma\bar{\delta}}
-2a^\beta_i\,\overline{a^\gamma_{\bar{j}}}\,\overline{a^\delta_{j}}\widetilde{R}^\alpha_{\beta\bar{\gamma}\bar{\delta}}.
\end{eqnarray}
On the other hand, by (\ref{2-1-16}), we have
\begin{eqnarray}\label{3-3}
\overline{\overline{a^\alpha_i}\,a^\beta_i\,\overline{a^\gamma_{\bar{j}}}\,\overline{a^\delta_{j}}\widetilde{R}^\alpha_{\beta\bar{\gamma}\bar{\delta}}}
&=&a^\alpha_i\,\overline{a^\beta_i}a^\gamma_{\bar{j}}\,a^\delta_{j}\widetilde{R}^\beta_{\alpha\delta\gamma}=\overline{a^\alpha_i}a^\beta_ia^\gamma_ja^\delta_{\bar{j}}\widetilde{R}^\alpha_{\beta\gamma\delta}.
\end{eqnarray}
The first identity follows from (\ref{3-1}), (\ref{3-2}) and (\ref{3-3}). Similar caculations do for the second identity.

\hfill$q.e.d$

The following Bochner formula belongs to V.Tosatti, which appeared in his work on Schwarz Lemma of pseudoholomorphic map between almost Hermitian manifolds.
\begin{Corollary}\label{Corollary4.2}\emph{(\cite{[T-07]})}
Let $f$ be a pseudoholomorphic map between almost Hermitian manifolds $(M,J,g)$ and $(\widetilde{M},\tilde{J}, \tilde{g})$, we have
\begin{eqnarray*}
\frac{1}{2}\Box \,e=|a^\alpha_{i,j}|^2+\overline{a^\alpha_i}a^\alpha_jR''_{i\bar{j}}-\overline{a^\alpha_i}a^\beta_ia^\gamma_j\overline{a^\delta_j}\widetilde{R}^\alpha_{\beta\gamma\bar{\delta}}.
\end{eqnarray*}
\end{Corollary}

\emph{Proof}. Since a pseudoholomorphic map is pluriharmonic, it follows from that $a^\alpha_{\bar{i}}=0$ and $e=e'$.

\hfill$q.e.d$





For the square of the length second fundamental form of pseudoholomorphic maps, we have
\begin{Proposition}\label{Proposition4.4}
Let $f$ be a pseudoholomorphic map between almost Hermitian manifolds $(M,J,g)$ and $(\widetilde{M},\tilde{J}, \tilde{g})$, and let $S$ be the half norm of the second fundamental form of $f$. Then
\begin{eqnarray*}
\frac{1}{2}\Box S&=&|a^\alpha_{i,jk}|^2+|a^\alpha_{i,j\bar{k}}|^2+ \overline{a^\alpha_{i,j}}\big(a^\alpha_{p,j}R''_{i\bar{p}}+a^\alpha_{i,p}R''_{j\bar{p}}+2a^\alpha_{p,k}R^p_{ij\bar{k}}-a^\beta_{i,j}a^\gamma_k\overline{a^\delta_k}\widetilde{R}^\alpha_{\beta\gamma\bar{\delta}}\big)\\
&{}&+2\emph{Re}\big(\overline{a^\alpha_{i,j}}a^\alpha_pR^p_{ij\bar{k},k}\big)-2\emph{Re}\big(\overline{a^\alpha_{i,j}}(a^\beta_ia^\gamma_j\overline{a^\delta_k}\widetilde{R}^\alpha_{\beta\gamma\bar{\delta}})_{,k}\big).
\end{eqnarray*}
\end{Proposition}

\emph{Proof}. Under the unitary frame field, we have $S=|a^\alpha_{i,j}|^2$, so
\begin{eqnarray}\label{4-4}
\frac{1}{2}\Box S&=&(a^\alpha_{i,j}\overline{a^\alpha_{i,j}})_{,k\bar{k}}\nonumber\\
&=&(a^\alpha_{i,jk}\overline{a^\alpha_{i,j}}+a^\alpha_{i,j}\overline{a^\alpha_{i,j\bar{k}}})_{,\bar{k}}\nonumber\\
&=& |a^\alpha_{i,jk}|^2+|a^\alpha_{i,j\bar{k}}|^2+(a^\alpha_{i,jk})_{,\bar{k}}\,\overline{a^\alpha_{i,j}}+a^\alpha_{i,j}\,(\overline{a^\alpha_{i,j\bar{k}}})_{,\bar{k}}.
\end{eqnarray}
By using the pseudoholomorphicity, the Ricci identities (\ref{2-3-80}) and (\ref{2-3-84}) give
\begin{equation}\label{4-5}
a^\alpha_{i,j\bar{k}}=a^\alpha_pR^p_{ij\bar{k}}-a^\beta_ia^\gamma_j\overline{a^\delta_k}\widetilde{R}^\alpha_{\beta\gamma\bar{\delta}},
\end{equation}
\begin{equation}\label{4-6}
(a^\alpha_{i,jk})_{,\bar{k}}=(a^\alpha_{i,j\bar{k}})_{,k}+a^\alpha_{p,j}R^p_{ik\bar{k}}+a^\alpha_{i,p}R^p_{jk\bar{k}}-a^\beta_{i,j}a^\gamma_k\overline{a^\delta_k}\widetilde{R}^\alpha_{\beta\gamma\bar{\delta}}.
\end{equation}
Then the half of the canonical Laplacian of $S$ is follows from (\ref{4-4}), (\ref{4-5}), (\ref{4-6}) and the curvature properties (2.15).

\hfill$q.e.d$

Appling the Bochner formulae above, we can get a $C^2$-estimate of pseudoholomorphic map. That is

\begin{Theorem}
Let $f$ be a pseudoholomorphic from closed almost Hermitian manifold $(M,J,g)$ into almost Hermitian manifold $(\widetilde{M},\tilde{J}, \tilde{g})$ with bounded geometry. If there is a constant $\Lambda$ such that $e(f)\leq\Lambda$, then
$$S\leq C,$$
where $C$ is a constant only depend on $\Lambda$, the curvatures and their covariant derivatives of $g$, $\tilde{g}$.
\end{Theorem}

\emph{Proof}. It follows from Corollary \ref{Corollary4.2} and Proposition \ref{Proposition4.4}, there exist constants $C_1$, $C_2$ and $C_3$ only depend on $\Lambda$, the correponding curvatures and their covariant derivatives such that 
\begin{equation*}
\frac{1}{2}\;\Box e(f)\geq S-C_1,\hspace{0.8cm}\frac{1}{2}\Box S\geq -C_2S-C_3,
\end{equation*}
which imply
\begin{equation*}
\frac{1}{2}\Box\big(S+(1+C_2)\,e(f)\big)\geq S-(C_1+C_1C_2+C_3).
\end{equation*}
By using the maximum principle to $S+(1+C_2)\,e(f)$ we can get the desired upper bound of $S$.

\hfill$q.e.d$

\section{Liouville type theorems}
We will use the properties of pseudoholomorphic maps and the Bochner formulae to get Liouville type results, which are similar to the theory of harmonic map. The main tools are the maximum principle and integral estimates.

\begin{Theorem}\label{Theorem5.1}
Let $f$ be a pseudoholomorphic map from closed almost Hermitian manifold $(M, J, g)$ into almost Hermitian manifold $(\widetilde{M},\tilde{J},\tilde{g})$, and $f(M)\subset V$. If $u$ is a strictly plurisubharmonic function on $V$, then $f$ is constant.
\end{Theorem}

\emph{Proof}. By using the pseudoholomorphicity, we have $a^\alpha_{\bar{i}}=0$ and $a^\alpha_{i,\bar{j}}=0$. So, the canonical Laplacian of $u\circ f$ is given by
\begin{eqnarray*}
\frac{1}{2}\Box (u\circ f)&=&(u\circ f)_{,k\bar{k}}\\
&=&(u_\alpha a^\alpha_k+u_{\bar{\alpha}}\, \overline{a^\alpha_{\bar{k}}})_{,\bar{k}}\\
&=&(u_\alpha a^\alpha_k)_{,\bar{k}}\\
&=&(u_{\alpha,\beta}\,a^\beta_{\bar{k}}+u_{\alpha,\bar{\beta}}\,a^{\bar{\beta}}_{\bar{k}})a^\alpha_k+u_\alpha a^\alpha_{k,\bar{k}}\\
&=&u_{\alpha,\bar{\beta}}\,a^\alpha_k \,\overline{a^\beta_k}.
\end{eqnarray*}
Notice that $u$ is strictly plurisubharmonic, so the maximum principle implies that $e(f)=0$, i.e., $f$ is a constant.

\hfill$q.e.d$

\begin{Theorem}
Let $f$ be a pseudoholomorphic map from closed almost Hermitian manifold $(M, J, g)$ into almost Hermitian manifold $(\widetilde{M},\tilde{J},\tilde{g})$
with $Ric''_M\geq 0$ and $HB_{\widetilde{M}}\leq 0$. Then $f$ is canonical totall geodesic. Furthermore, $f$ is constant if $Ric''_M$ is positive at a point in $M$.
\end{Theorem}

\emph{Proof}. Recall the Bochner formule obtained in Corollary \ref{Corollary4.2}, that is
\begin{equation}\label{5-1}
\frac{1}{2}\Box \,e=|a^\alpha_{i,j}|^2+\overline{a^\alpha_i}a^\alpha_jR''_{i\bar{j}}-\overline{a^\alpha_i}a^\beta_ia^\gamma_j\overline{a^\delta_j}\widetilde{R}^\alpha_{\beta\gamma\bar{\delta}}.
\end{equation}
The curvature conditons imply that $\Box e\geq 0$, so the maximum principle implies $e$ is a constant. At this time, all the terms in the right hand side of (\ref{5-1}) are nonnegative and their summation is equal to zero, so we obtain $|a^\alpha_{i,j}|^2=0$, i.e., $\nabla df=0$. If $Ric''_M$ is positive at a piont, the second term in the right hand side of (\ref{5-1}) implies $e=0$.

\hfill$q.e.d$ 

We can also obtain a similar result when the holomorphic bisectional curvature of target manifold might be positive.  

\begin{Theorem}\label{Theorem5.3}
Let $(M, J, g)$, $(\widetilde{M},\tilde{J},\tilde{g})$ be two  almost Hermitian manifolds with $M$ is closed, and $Ric''_M\geq A$, $HB_{\widetilde{M}}\leq B$, where $A$, $B$ are positive constants. If the pseudoholomorphic map $f:M\longrightarrow\widetilde{M}$ satisfies $e(f)\leq \frac{A}{B}$, then $f$ is canonical totally geodesic and $e(f)=0$ or $e(f)=\frac{A}{B}$. Furthermore, $f$ is constant if $e(f)< \frac{A}{B}$.
\end{Theorem}

\emph{Proof}. By the curvture conditions, the Bochner formule (\ref{5-1}) implies
\begin{equation*}
\frac{1}{2}\Box\,e(f)\geq |a^\alpha_{i,j}|^2+\big(A-B\,e(f)\big)e(f).
\end{equation*}
So, if $e(f)\leq \frac{A}{B}$, we have $\Box\,e(f)\geq 0$, then the maximum principle implies $e(f)$ is a constant, $|a^\alpha_{i,j}|^2=0$ and $\big(A-B\,e(f)\big)e(f)=0$. Thus, we have $e(f)=0$ if  $e(f)< \frac{A}{B}$, which means that $f$ is constant.

\hfill$q.e.d$

For the domain manifold is complete and non-compact, we can obtain the following Liouville type result provied the energy is finite.
\begin{Theorem}
Let $(M, J, g)$ be a complete, non-compact almost Hermitian manifold with infinite volume and $Ric''_M\geq 0$, let $(\widetilde{M},\tilde{J},\tilde{g})$ be an almost Hermitian manifold with $HB_{\widetilde{M}}\leq 0$. If the pseudoholomorphic $f:M\longrightarrow\widetilde{M}$ has finite energy and 
$$\int_{M}\langle \nabla e(f), X_L\rangle\,dV_g=0,$$
then $f$ is constant.
\end{Theorem}

\emph{Proof}. The Bochner formule (\ref{5-1}) and Cauchy-Schwarz inequality imply
\begin{equation}\label{5-2}
\Box\,e(f)\geq 2S, \hspace{1cm} |\nabla\,e(f)|^2\leq 2\,e(f)\,S,
\end{equation}
respectively. For any $\epsilon >0$, set $u_\epsilon:=\sqrt{e(f)+\epsilon}$, then we have
\begin{eqnarray*}
\frac{1}{2}\Box u_\epsilon&=&(u_\epsilon)_{,k\bar{k}}\\
&=& \frac{1}{4}u_\epsilon^{-1}\big(\Box\, e(f)-\frac{|\nabla e(f)|^2}{2(e(f)+\epsilon)}\big).
\end{eqnarray*}
Together with (\ref{5-2}), this implies
\begin{equation}\label{5-3}
\Box u_\epsilon\geq \frac{1}{2}u_\epsilon^{-1}S\big(1-\frac{e(f)}{2(e(f)+\epsilon)}\big)\geq0.
\end{equation}

Let $x_0\in M$ be a fixed point. We denote the geodesic balls centered at $x_0$ with radii $R$, $2R$ by $B_R$, $B_{2R}$, respectively. Choosing a cut-off function $\eta$ such that $0\leq \eta \leq1$, $\eta(x)=1$ for $x\in B_R$, $\eta(x)=0$ for $x\in M\setminus B_{2R}$. Moreover, we further assume that $|\nabla \eta|\leq \frac{C}{R}$, where $C$ is a positive constant. By using (2.24), the inequality (\ref{5-3}) gives
\begin{eqnarray}\label{5-4}
0&\leq& \int_{B_{2R}}\eta^2\, u_\epsilon\,\Box u_\epsilon\,dV_g\nonumber\\
&=&\int_{B_{2R}}\eta^2\, u_\epsilon\,\big(\Delta u_\epsilon+2\langle  \nabla u_\epsilon, X_L\rangle\big)\,dV_g\nonumber\\
&=&-2\int_{B_{2R}}\eta\,u_\epsilon\langle \nabla\eta,\nabla u_\epsilon\rangle\, dV_g-\int_{B_{2R}}\eta^2|\nabla u_\epsilon|^2+\int_{B_{2R}}\eta^2\langle \nabla e(f),X_L\rangle\,dV_g\nonumber\\
&\leq& 2\Big(\int_{B_{2R}\setminus B_R}\eta^2|\nabla u_\epsilon|^2\,dV_g\Big)^{1/2}\Big(\int_{B_{2R}\setminus B_R}u_\epsilon^2|\nabla \eta|^2\,dV_g\Big)^{1/2}
-\int_{B_{2R}\setminus B_R} \eta^2|\nabla u_\epsilon|^2\,dV_g\nonumber\\
&{}&-\int_{B_R} |\nabla u_\epsilon|^2\,dV_g+\int_{B_{2R}}\eta^2\langle \nabla e(f),X_L\rangle\,dV_g.
\end{eqnarray}
If we view (\ref{5-4}) as a quadratic inequality of the term 
$\Big(\int_{B_{2R}\setminus B_R} \eta^2|\nabla u_\epsilon|^2\,dV_g\Big)^{1/2}$, we get
\begin{eqnarray}\label{5-5}
\int_{B_R} |\nabla u_\epsilon|^2\,dV_g-\int_{B_{2R}}\eta^2\langle \nabla e(f),X_L\rangle\,dV_g&\leq& \int_{B_{2R}\setminus B_R}u_\epsilon^2|\nabla \eta|^2\,dV_g\nonumber\\
&\leq& \frac{C^2}{R^2}\int_{B_{2R}} u^2_\epsilon\;dV_g.
\end{eqnarray}
Notice that $|\nabla u_\epsilon|^2=\frac{|\nabla e(f)|^2}{4(e(f)+\epsilon)}$, the inequality (\ref{5-5}) implies
\begin{eqnarray}\label{5-6}
\int_{\Sigma \cap B_R} \frac{|\nabla e(f)|^2}{4(e(f)+\epsilon)}\,dV_g-\int_{B_{2R}}\eta^2\langle \nabla e(f),X_L\rangle\,dV_g
&\leq& \frac{C^2}{R^2}\int_{B_{2R}} (e(f)+\epsilon)\;dV_g,
\end{eqnarray}
where $\Sigma:=\{x\in M\,|\,e(f)(x)\neq 0\}.$ Letting $\epsilon\rightarrow 0$, $R\rightarrow +\infty$ and using $f$ has finite energy, we obtain
\begin{eqnarray}\label{5-6}
\int_{\Sigma} \frac{|\nabla e(f)|^2}{4e(f)}\,dV_g
\leq 0.
\end{eqnarray}
This tells us that $e(f)$ is a constant. Hence, $e(f)=0$ from the facts that $f$ has finite energy and $(M,g)$ has infinite volume.

\hfill$q.e.d$

\emph{Remark.} The condition that $(M, g)$ has infinite volume can be replaced by a curvature condition. The integration condition is always stisfied when $(M,J,g)$ is semi-K\"ahlerian, quasi-K\"ahlerian, or almost-K\"ahlerian, or nearly-K\"ahlerian by the fact that $X_L=0$.

\section{Pseudoholomorphicity of pluriharmonic maps}
We will investigate the pseudoholomorphicity of pluriharmonic maps by using the Bochner formulae and maximum principle. Two cases with respect to the topology of domain manifold will be considered.

\begin{Theorem}
Let$(M, J, g)$ be a closed almost Hermitian manifold with $Ric''_M\geq A$, and let  $(\widetilde{M},\tilde{J},\tilde{g})$ be an almost Hermitian manifold with $HB_{\widetilde{M}}\leq B$ and the $(2,0)$-part of the curvature is bounded by $C$. If a pluriharmonic map $f:M\longrightarrow\widetilde{M}$ satisfies $(B+2C)\,e(f)< A$,
then $f$ is pseudoholomorphic.
\end{Theorem}

\emph{Proof}. Under the curvature conditions, the second Bochner formule in Proposition 3.1 gives
\begin{eqnarray}\label{6-1}
\frac{1}{2}\Box e''(f)&\geq& |a^\alpha_{\bar{i},\bar{j}}|^2+e''(f)\big(A-4C\,\sqrt{e''(f)}\sqrt{e'(f)}-B\,e(f)\big)\nonumber\\
&\geq&  |a^\alpha_{\bar{i},\bar{j}}|^2+e''(f)\Big(A-2C\,e(f)-B\,e(f)\Big).
\end{eqnarray}
By using the maximum principle, the differential inequality (\ref{6-1}) implies 
$$e''(f)\big(A-2C\,e(f)-B\,e(f)\big)=0.$$
Thus, we have $e''(f)=0$ by the fact that $(B+2C)\,e(f)< A$.

\hfill$q.e.d$

To proceed, we need the following maximum principle (Proposition 4.1 in \cite{[T-07]}) due to V.Tosatti.

\begin{Lemma}\label{Lemma6.2}\emph{(\cite{[T-07]})}
Let $(M,J,g)$ be a complete almost Hermitian manifold with second Ricci curvature bounded below, with torsion and $(2,0)$-part of the curvature bounded. Let $u$ be a real function that is bounded from below. Then ginven any $\epsilon>0$ there exists a point $x_\epsilon\in M$ such that
\begin{equation*}
\liminf\limits_{\epsilon\rightarrow 0}u(x_\epsilon)=\inf\limits_{M}u,\hspace{0.6cm}
|\nabla u|(x_\epsilon)\leq \epsilon,\hspace{0.6cm}
\Box\, u(x_\epsilon)\geq -\epsilon.
\end{equation*}
\end{Lemma}

\begin{Theorem}\label{Theorem6.3}
Let $(M,J,g)$ be a complete almost Hermitian manifold with $Ric''_M\geq A$ and with torsion and $(2,0)$-part of the curvature bounded by $C$. Let  $(\widetilde{M},\tilde{J},\tilde{g})$ be an almost Hermitian manifold with $B_1\leq
HB_{\widetilde{M}}\leq B_2$. If a pluriharmonic map $f:M\longrightarrow \widetilde{M}$ satisfies $(B_2+2C)e''(f)\leq A+(B_1-2C)e'(f)-\delta$ for a positive constant $\delta>0$, then $f$ is pseudoholomorphic.
\end{Theorem}

\emph{Proof}. For any $\lambda>0$, set $u_\lambda:=\big(e''(f)+\lambda\big)^{-1/2}$, through direct calculation, we have
\begin{equation}\label{6-2}
\Box\,u_\lambda=3\,u_\lambda^{-1}\,|\nabla u_\lambda|^2-\frac{1}{2}u_\lambda^{3}\,\Box e''(f).
\end{equation}
On the other hand, by using the second Bochner formula in Proposition 4.1, we obtain
\begin{eqnarray}\label{6-3}
\frac{1}{2}\Box\,e''(f)&\geq& A\,e''(f)-4C\,\big(e'(f)\big)^{1/2}\big(e''(f)\big)^{3/2}+B_1\,e'(f)\,e''(f)-B_2\big(e''(f)\big)^2\nonumber\\
&=& e''(f)\Big[A-4C\,\big(e'(f)\big)^{1/2}\big(e''(f)\big)^{1/2}+B_1\,e'(f)-B_2\,e''(f)\Big]\nonumber\\
&\geq& e''(f)\big[A+(B_1-2C)\,e'(f)-(B_2+2C)\,e''(f)\big]\nonumber\\
&\geq& \delta\,e''(f).
\end{eqnarray}
Applying the Lemma 5.2 to $u_\lambda$, for any $\epsilon>0$, there exisits a point $x_\epsilon\in M$ such that
\begin{equation*}
u_\lambda(x_\epsilon)\leq\inf\limits_{M}u_\lambda+\epsilon,\hspace{0.6cm}
|\nabla u_\lambda|(x_\epsilon)\leq \epsilon,\hspace{0.6cm}
\Box\, u_\lambda(x_\epsilon)\geq -\epsilon.
\end{equation*}
So, at point $x_\epsilon$, the identity (\ref{6-2}) and the inequality (\ref{6-3}) imply
\begin{eqnarray}\label{6-4}
\frac{\delta\,e''(f)}{(e''(f)+\lambda)^2}=\delta\,e''(f)\,u_\lambda^{4}&\leq& 3\epsilon^2+\epsilon\,u_\lambda\nonumber\\
&\leq&3\epsilon^2+\epsilon\big(\inf\limits_{M}u_\lambda+\epsilon\big).
\end{eqnarray}
Letting $\epsilon\rightarrow 0$, then $u_\lambda(x_\epsilon)\rightarrow \inf\limits_{M}u_\lambda$ and hence $e''(f)(x_\epsilon)\rightarrow\sup\limits_M e''(f)$. Thus, the inequality (\ref{6-4}) implies that $\sup\limits_M e''(f)\leq 0$, i.e., $f$ is pseudoholomorphic.

\hfill$q.e.d$

\section{Pseudoholomorphic isometric immersion}
The purpose of this section is to derive Simons integral inequality, to give bounds of the norm of parallel canonical second fundamental form, of the   pseudoholomorphic isometry between the almost Hermitian manifolds. The idear inspired from H.W.Xu's work (\cite{[X-95]}) on closed minimal submanifold in pinched Riemannian manifolds. We will adopt the following range of indices:
$$1\leq i,j,k,\ell\ldots\leq m;\hspace{0.3cm}m+1\leq \lambda,\mu,\nu,\sigma,\ldots\leq n; \hspace{0.3cm}1\leq \alpha,\beta,\gamma,\delta\ldots\leq n.$$

We first recall the structure equations of almost Hermitian manifold $(\widetilde{M},\tilde{J},\tilde{g})$. Locally, by choosing a unitary frame field $\{\tilde{e}_\alpha\}$ with dual $\{\tilde{\theta}^\alpha\}$, the first and second structure equations are given by
\begin{equation}\label{7-1}
d\tilde{\theta}^\alpha=-\tilde{\theta}^\alpha_\beta\wedge\tilde{\theta}^\beta+\widetilde{\Theta}^\alpha,\hspace{0.3cm}\tilde{\theta}^\alpha_\beta+\overline{\tilde{\theta}^\beta_\alpha}=0,
\end{equation}
\begin{equation}\label{7-2}
d\tilde{\theta}^\alpha_{\beta}=-\tilde{\theta}^\alpha_\gamma\wedge\tilde{\theta}^\gamma_\beta+\widetilde{\Omega}^\alpha_\beta,\hspace{0.3cm}\widetilde{\Omega}^\alpha_\beta+\overline{\widetilde{\Omega}^\beta_\alpha}=0,
\end{equation}
where $\tilde{\theta}^\alpha_\beta$, $\widetilde{\Theta}^\alpha$, $\widetilde{\Omega}^\alpha_\beta$ are connection 1-forms, torsion and curvature forms respectively.

Let $f$ be a pseudoholomorphic isometry from $(M,J,g)$ into $(\widetilde{M},\tilde{J},\tilde{g})$ with $\mbox{dim}M=2m$ and  $\mbox{dim}\widetilde{M}=2n$. We choose an adapted unitary frame field $\{\tilde{e}_\alpha\}$ along $f$, which means that $\{\tilde{e}_i\}$ are tangent to $f(M)$ and $\{\tilde{e}_\lambda\}$ are normal to $f(M)$. So, we have
\begin{equation}\label{7-3}
f^*\tilde{\theta}^i=\theta^i,\hspace{1cm}f^*\tilde{\theta}^\lambda=0,
\end{equation}
i.e., $a^\alpha_i=\delta_{\alpha\,i}$. Then the pull-back of (\ref{7-1}) with $\alpha=i$ gives
\begin{equation}\label{7-4}
\tilde{\theta}^i_j=\theta^i_j,\hspace{1cm}\widetilde{\Theta}^i=\Theta^i,
\end{equation}
where we have omitted the pull-back $f^*$ on the forms $\tilde{\theta}^i_j$, $\widetilde{\Theta}^i$. Substituting $\tilde{\theta}^i_j=\theta^i_j$ and $a^\alpha_i=\delta_{\alpha i}$ into (\ref{2-3-73}) with $\alpha=k$, $\lambda$ respectively, we have
\begin{equation}\label{7-5}
a^k_{i,j}=0,\hspace{1cm}\tilde{\theta}^\lambda_i=a^\lambda_{i,j}\,\theta^j.
\end{equation}
Taking $(\alpha,\beta)=(i,j)$ in (\ref{7-2}), together with the first identities in(\ref{7-4}), we obtain the Gauss equation
\begin{equation*}
\Omega^i_j+(-\overline{\tilde{\theta}^\lambda_i})\wedge\tilde{\theta}^\lambda_j=\widetilde{\Omega}^i_j,
\end{equation*}
which implies
\begin{equation}\label{7-6}
R^i_{jk\bar{\ell}}+a^\lambda_{j,k}\overline{a^\lambda_{i,\ell}}=\widetilde{R}^i_{jk\bar{\ell}},
\end{equation}
and
\begin{equation}\label{7-7}
R^i_{jk\ell}=\widetilde{R}^i_{jk\ell},\hspace{0.6cm} R^i_{j\bar{k}\;\bar{\ell}}=\widetilde{R}^i_{j\bar{k}\;\bar{\ell}}.
\end{equation}
Similarly, taking $(\alpha,\beta)=(\lambda,\mu)$ in (\ref{7-2}), we get the Ricci equation
\begin{equation}\label{7-8}
\Omega^\lambda_\mu+\tilde{\theta}^\lambda_i\wedge(-\overline{\tilde{\theta}^\mu_i})=\widetilde{\Omega}^\lambda_\mu,
\end{equation}
where $\Omega^\lambda_\mu$ is the normal curvature forms defined by $\Omega^\lambda_\mu:=d\tilde{\theta}^\lambda_\mu+\tilde{\theta}^\lambda_\nu\wedge\tilde{\theta}^\nu_\mu$. If we set
\begin{equation*}
\Omega^\lambda_\mu:=R^\lambda_{\mu jk}\,\theta^j\wedge\theta^k +R^\lambda_{\mu j\bar{k}}\,\theta^j\wedge\overline{\theta^k} +R^\lambda_{\mu \bar{j}\,\bar{k}}\,\overline{\theta^j}\wedge\overline{\theta^k}, 
\end{equation*}
from (\ref{7-8}), we have
\begin{equation}\label{7-9}
R^\lambda_{\mu j\bar{k}}=a^\lambda_{i,j}\overline{a^\mu_{i,k}}+\widetilde{R}^\lambda_{\mu j\bar{k}},
\end{equation}
\begin{equation}\label{7-10}
R^\lambda_{\mu j k}=\widetilde{R}^\lambda_{\mu j k},\hspace{0.6cm}R^\lambda_{\mu \bar{j}\, \bar{k}}=\widetilde{R}^\lambda_{\mu \bar{j}\, \bar{k}}.
\end{equation}
We call $R^\perp:=\sum\limits_{i,\lambda}R^\lambda_{\lambda i i}$ the \emph{normal scalar curvature} of $M$ in $\widetilde{M}$.
The Codazzi equation is follows from (\ref{2-3-79})-(\ref{2-3-81}) or by taking $(\alpha, \beta)=(\lambda,i)$ in (\ref{7-2}), that is
\begin{equation}\label{7-11}
a^\lambda_{i,j\bar{k}}=-\widetilde{R}^\lambda_{ij\bar{k}}, \hspace{0.5cm}
a^\lambda_{i,jk}-a^\lambda_{i,kj}=2a^\lambda_{i,p}L^p_{jk}-2\widetilde{R}^\lambda_{ijk},\hspace{0.5cm}
a^\lambda_{i,p}N^p_{\bar{j}\,\bar{k}}=\widetilde{R}^\lambda_{i\bar{j}\,\bar{k}}.
\end{equation}

There are some basic properties of pseudoholomorphic isometry.

\begin{Proposition}
Let $(M, J, g)$, $(\widetilde{M},\tilde{J},\tilde{g})$ be two almost Hermitian manifolds with $HB_{\widetilde{M}}\leq B$. If $f:M \rightarrow \widetilde{M}$ is a pseudoholomorphic isometry, then we have

\emph{(1)} $Ric'_M-mBg$, $Ric''_M-mBg$ is negative semi-definite.

\emph{(2)} $R\leq m^2 B$.

\emph{(3)} ${HS}_M\leq B$.
\end{Proposition}

\emph{Proof}. For any $X=X^i\,e_i\in T'M$, by using the Gauss equation (\ref{7-6}), we have
\begin{eqnarray*}
R'_{j\bar{k}}X^j\overline{X^k}-mB|X|^2&\leq&\sum\limits_{i}R^i_{ij\bar{k}}X^j\overline{X^k}-\sum\limits_{i}\widetilde{R}^i_{ij\bar{k}}X^j\overline{X^k}\\
&=&-\sum\limits_{i,\lambda}a^\lambda_{i,j}\overline{a^\lambda_{i,k}}X^j\overline{X^k}\leq 0.
\end{eqnarray*} 
Similar proof for $Ric''_M-mBg$. The second and third inequalities are also follwow from Gauss equation (\ref{7-6}) and  $HB_{\widetilde{M}}\leq B$.

\hfill$q.e.d$

\emph{Remark}. The Proposition 6.1 provides various curvature criterion to detect that whether an almost Hermitian manifold can be pseudoholomorphically immersed into another one.  

\begin{Proposition}
Let $(M, J, g)$, $(\widetilde{M},\tilde{J},\tilde{g})$ be two almost Hermitian manifolds with $HB_{\widetilde{M}}\leq B$. If $f:M \rightarrow \widetilde{M}$ is a pseudoholomorphic isometry, then $f$ is canonical totally geodesic if $f$ satisfies one of the following conditions\emph{:}

\emph{(1)} $Ric'_M=mBg$, or $Ric''_M=mBg$.

\emph{(2)} $R= m^2 B$.

\emph{(3)} ${HS}_M=B$.
\end{Proposition}

\emph{Proof}. By the Gauss equation (\ref{7-6}), any one of the conditions implies that the canonical second fundamental form $a^\alpha_{i,j}=0$, i.e., 
$\nabla df=0$.

\hfill$q.e.d$

For the scalar curvature and normal scalar curvature, we have
 
\begin{Proposition}
Let $(M, J, g)$, $(\widetilde{M},\tilde{J},\tilde{g})$ be two almost Hermitian manifolds with $B_1\leq HB_{\widetilde{M}}\leq B_2$. For a pseudoholomorphic isometry $f:M \rightarrow \widetilde{M}$, we have

$$mnB_1\leq R+R^\perp\leq mnB_2.$$
\end{Proposition}

\emph{Proof}. By the Gauss equation (\ref{7-6}) and Ricci equation (\ref{7-9}), we have
$$R+R^\perp=\sum\limits_{i,j}\widetilde{R}^i_{ij\bar{j}}+\sum\limits_{i,\lambda}\widetilde{R}^\lambda_{\lambda i\bar{i}}.$$ 
So, the inequalities follow from the bounds of holomorphic bisectional curvature of $\widetilde{M}$.

\hfill$q.e.d$

We will give an upper bound of the curvature tensor in terms of the bound of holomorphic bisectional curvature. 

\begin{Lemma}\label{Lemma6.4}
Let $(M, J, g)$ be an almost Hermitian manifold with $a \leq HB_M\leq b$. Under a unitary frame field, the components of $(1,1)$-part of the curvature denoted by $R^i_{jk\bar{\ell}}$, then for fixed indices $i,j,k,\ell$, we have

\emph{(1)} $|R^j_{ik\bar{\ell}}X^i\overline{X^j}|^2\leq 2 (b-a)^2|X|^4$ for $k\neq\ell$ and any $X\in T'M$.

\emph{(2)} $|R^i_{jk\bar{\ell}}X^k\overline{X^\ell}|^2\leq 2 (b-a)^2|X|^4$ for $i\neq j$ and any $X\in T'M$.

\emph{(3)} $|R^j_{ik\bar{\ell}}|^2 \leq 4(1+\sqrt{2})^2(b-a)^2$ for $i\neq j$ and $k\neq\ell$.
\end{Lemma}

\emph{Proof}. For any $X=X^i\,e_i$, $Y=Y^j\,e_j$, $Z=Z^k\,e_k\in T'M$, we have
\begin{eqnarray}\label{7-12}
R^j_{ik\bar{\ell}}\,X^i\overline{X^j}Y^k\overline{Z^\ell}+R^j_{ik\bar{\ell}}\,X^i\overline{X^j}Z^k\overline{Y^\ell}&=&R^j_{ik\bar{\ell}}\,X^i\overline{X^j}(Y^k+Z^k)(\overline{Y^\ell}+\overline{Z^\ell})\hspace{2cm}\nonumber\\
&{}&-R^j_{ik\bar{\ell}}\,X^i\overline{X^j}Y^k\overline{Y^\ell}-R^j_{ik\bar{\ell}}\,X^i\overline{X^j}Z^k\overline{Z^\ell}.
\end{eqnarray}
Taking $X=X^i\,e_i$, $Y=e_k$, $Z=e_\ell$ with $k\neq \ell$ in (\ref{7-12}) and together $\overline{R^i_{jk\bar{\ell}}}=R^j_{i\ell\,\bar{k}}$, we obtain
\begin{eqnarray}\label{7-13}
2\,\mbox{Re}\big(R^j_{ik\bar{\ell}}X^i\overline{X^j}\big)&=&R^j_{ik\bar{\ell}}X^i\overline{X^j}+R^j_{i\ell\,\bar{k}}X^i\overline{X^j}\hspace{7cm}\nonumber\\
&=& R^j_{ip\,\bar{q}}X^i\overline{X^j}(\delta_{pk}+\delta_{p\ell})(\delta_{qk}+\delta_{q\ell})-R^j_{ik\,\bar{k}}X^i\overline{X^j}-R^j_{i\ell\,\bar{\ell}}X^i\overline{X^j},
\end{eqnarray}
which implies
\begin{equation}\label{7-14}
-(b-a)|X|^2\leq \mbox{Re}(R^j_{ik\bar{\ell}}X^i\overline{X^j})\leq (b-a)|X|^2.
\end{equation}
By taking $X=X^ie_i$, $Y=-\sqrt{-1}e_k$, $Z=e_\ell$ with $k\neq \ell$ in (\ref{7-12}), we can obtain
\begin{equation}\label{7-15}
-(b-a)|X|^2\leq \mbox{Im}(R^j_{ik\bar{\ell}}X^i\overline{X^j})\leq (b-a)|X|^2.
\end{equation}
So, inequlities (\ref{7-14}) and (\ref{7-15}) give the desired upper bound of $|R^j_{ik\bar{\ell}}X^i\overline{X^j}|^2$, and similar proof for $|R^i_{jk\bar{\ell}}X^k\overline{X^\ell}|^2$.

On the other hand, for any $X=X^i\,e_i$, $Y=Y^j\,e_j$, $Z=Z^k\,e_k$, $W=W^\ell\,e_\ell\in T'M$, we have
\begin{eqnarray}\label{7-16}
R^j_{ik\bar{\ell}}\,X^i\overline{Y^j}Z^k\overline{W^\ell}+R^j_{ik\bar{\ell}}\,Y^i\overline{X^j}Z^k\overline{W^\ell}&=&R^j_{ik\bar{\ell}}(X^i+Y^i)(\overline{X^j}+\overline{Y^j})Z^k\overline{W^\ell}\hspace{2cm}\nonumber\\
&{}&-R^j_{ik\bar{\ell}}\,X^i\overline{X^j}Z^k\overline{W^\ell}-R^j_{ik\bar{\ell}}\,Y^i\overline{Y^j}Z^k\overline{W^\ell}.
\end{eqnarray}
Taking $X=e_i$, $Y=e_j$, $Z=e_k$, $W=e_\ell$ with  $i\neq j$ and $k\neq \ell$ in (\ref{7-16}), we obtain
\begin{eqnarray}\label{7-17}
R^j_{ik\bar{\ell}}+R^i_{jk\bar{\ell}}=R^q_{pk\bar{\ell}}\,(\delta_{pi}+\delta_{pj})(\delta_{qi}+\delta_{qj})-R^i_{ik\bar{\ell}}-R^j_{jk\bar{\ell}},
\end{eqnarray}
and by taking $X=\sqrt{-1}e_i$, $Y=e_j$, $Z=e_k$, $W=\sqrt{-1}e_\ell$ with  $i\neq j$ and $k\neq \ell$ in (\ref{7-16}), we obtain
\begin{eqnarray}\label{7-18}
R^j_{ik\bar{\ell}}-R^i_{jk\bar{\ell}}=-\sqrt{-1}R^q_{pk\bar{\ell}}\,(\sqrt{-1}\delta_{pi}+\delta_{pj})(-\sqrt{-1}\delta_{qi}+\delta_{qj})+\sqrt{-1}R^i_{ik\bar{\ell}}+\sqrt{-1}R^j_{jk\bar{\ell}}.\hspace{0.3cm}
\end{eqnarray}
Thus, from (\ref{7-17}) and (\ref{7-18}), we get
\begin{eqnarray}\label{7-19}
2R^j_{ik\bar{\ell}}&=&R^q_{pk\bar{\ell}}\,(\delta_{pi}+\delta_{pj})(\delta_{qi}+\delta_{qj})-\sqrt{-1}R^q_{pk\bar{\ell}}\,(\sqrt{-1}\delta_{pi}+\delta_{pj})(-\sqrt{-1}\delta_{qi}+\delta_{qj})\nonumber\\
&{}&-(1-\sqrt{-1})R^i_{ik\bar{\ell}}-(1-\sqrt{-1})R^j_{jk\bar{\ell}}.
\end{eqnarray}
So, by using the estimates (1) and (2) have been obtained, we get
\begin{eqnarray}
2|R^j_{ik\bar{\ell}}| &\leq& |R^q_{pk\bar{\ell}}\,(\delta_{pi}+\delta_{pj})(\delta_{qi}+\delta_{qj})|+|R^q_{pk\bar{\ell}}\,(\sqrt{-1}\delta_{pi}+\delta_{pj})(-\sqrt{-1}\delta_{qi}+\delta_{qj})|\nonumber\\
&{}& +\sqrt{2}\,|R^i_{ik\bar{\ell}}|+\sqrt{2}\,|R^j_{jk\bar{\ell}}|\nonumber\\
&\leq& \sqrt{2}(b-a)|e_i+e_j|^2+ \sqrt{2}(b-a)|\sqrt{-1}e_i+e_j|^2 +4(b-a)\nonumber\\
&=&4(1+\sqrt{2})(b-a),
\end{eqnarray}
which gives the desired estimate of $|R^j_{ik\ell}|^2$ with $i\neq j$ and $k\neq\ell$.

\hfill$q.e.d$

We wish to give a lower bound of $\frac{1}{2}\Box S$ by using the bounds of holomorphic bisectional curvature of target manifold $(\widetilde{M},\tilde{J},\tilde{g})$. In the sequal, we assume that the holomorphic bisectional curvature pointwisely satisfies $a(x)\leq HB_{\widetilde{M}}\leq b(x)$ for $x\in \widetilde{M}$. According to the expression of $\frac{1}{2}\Box S$, for convenience, we define three terms as follows:
\begin{equation*}
\mbox{(I)}:= 2\mbox{Re}\big(\overline{a^\alpha_{i,j}}a^\alpha_pR^p_{ij\bar{k},k}\big),
\end{equation*}
\begin{equation*}
\mbox{(II)}:= |a^\alpha_{i,j\bar{k}}|^2-2\mbox{Re}\big(\overline{a^\alpha_{i,j}}(a^\beta_ia^\gamma_j\overline{a^\delta_k}\widetilde{R}^\alpha_{\beta\gamma\bar{\delta}})_{,k}\big),
\end{equation*}
\begin{equation*}
\mbox{(III)}:= \overline{a^\alpha_{i,j}}\big(a^\alpha_{p,j}R''_{i\bar{p}}+a^\alpha_{i,p}R''_{j\bar{p}}+2a^\alpha_{p,k}R^p_{ij\bar{k}}-a^\beta_{i,j}a^\gamma_k\overline{a^\delta_k}\widetilde{R}^\alpha_{\beta\gamma\bar{\delta}}\big).
\end{equation*}

\begin{Lemma}\label{Lemma6.5}
For the terms \emph{(I)} and \emph{(II)}, we have

\emph{(1)} $\emph{(I)}=0$.

\emph{(2)} $\emph{(II)}\geq -\mbox{\emph{div}}^c(X)- P(m,n)(b-a)^2-m^2\,\emph{max}\{a^2,b^2\}$, where $X$ defined in \emph{(\ref{6.22})} and $P(m,n)$ is a polynomial with respect to $m,n$ defined in \emph{(\ref{6.25})}.
\end{Lemma}

\emph{Proof}. By choosing an adapted unitary frame field, we have $a^k_{i,j}=0$ and $a^\lambda_i=0$, which imply $\mbox{(I)}=0$. Notice that 
\begin{eqnarray}\label{6.21}
2\mbox{Re}\big(\overline{a^\alpha_{i,j}}(a^\beta_ia^\gamma_j\overline{a^\delta_k}\widetilde{R}^\alpha_{\beta\gamma\bar{\delta}})_{,k}\big)
&=& \overline{a^\alpha_{i,j}}(a^\beta_ia^\gamma_j\overline{a^\delta_k}\widetilde{R}^\alpha_{\beta\gamma\bar{\delta}})_{,k}+a^\alpha_{i,j}(\overline{a^\beta_ia^\gamma_j\overline{a^\delta_k}\widetilde{R}^\alpha_{\beta\gamma\bar{\delta}}})_{,\bar{k}}\nonumber\\
&=&X^k_{,k}
+\overline{X^k_{,k}}-\overline{a^\alpha_{i,j\bar{k}}}\widetilde{R}^\alpha_{ij\bar{k}}-a^\alpha_{i,j\bar{k}}\overline{\widetilde{R}^\alpha_{ij\bar{k}}}\nonumber\\
&=&\mbox{div}^c(X)-\overline{a^\alpha_{i,j\bar{k}}}\widetilde{R}^\alpha_{ij\bar{k}}-a^\alpha_{i,j\bar{k}}\overline{\widetilde{R}^\alpha_{ij\bar{k}}}
\end{eqnarray}
where 
\begin{equation}\label{6.22}
X=X^k\,e_k+\overline{X^k}\,\overline{e_k},\hspace{0.6cm}X^k=\overline{a^\alpha_{i,j}}a^\beta_ia^\gamma_j\overline{a^\delta_k}\widetilde{R}^\alpha_{\beta\gamma\bar{\delta}}.
\end{equation}
 Thus, by (\ref{6.21}) and Cauchy-Schwarz inequality, we have
\begin{eqnarray}\label{6.23}
\mbox{(II)}&=& |a^\alpha_{i,j\bar{k}}|^2-\mbox{div}^c(X)+\overline{a^\alpha_{i,j\bar{k}}}\widetilde{R}^\alpha_{ij\bar{k}}+a^\alpha_{i,j\bar{k}}\overline{\widetilde{R}^\alpha_{ij\bar{k}}}\nonumber\\
&\geq& |a^\alpha_{i,j\bar{k}}|^2-\mbox{div}^c(X)-|a^\alpha_{i,j\bar{k}}|^2-|\widetilde{R}^\alpha_{i,j\bar{k}}|^2\nonumber\\
&=&-\mbox{div}^c(X)-|\widetilde{R}^\alpha_{ij\bar{k}}|^2.
\end{eqnarray}
On the other hand, by using Lemma 6.4 repeatedly, we have
\begin{eqnarray}\label{6.24}
\sum\limits_{i,j,k,\alpha}|\widetilde{R}^\alpha_{ij\bar{k}}|^2&=&\sum\limits_{i,j\neq k}|\widetilde{R}^i_{ij\bar{k}}|^2+\sum\limits_{i\neq\ell,j}|\widetilde{R}^\ell_{ij\bar{j}}|^2+\sum\limits_{i,j}|\widetilde{R}^i_{ij\bar{j}}|^2+\sum\limits_{i\neq\ell,j\neq k}|\widetilde{R}^\ell_{ij\bar{k}}|^2\nonumber\\
&{}&+\sum\limits_{i,j,\lambda}|\widetilde{R}^\lambda_{ij\bar{j}}|^2+\sum\limits_{i,j\neq k,\lambda}|\widetilde{R}^\lambda_{ij\bar{k}}|^2\nonumber\\
&\leq& 2m^2(m-1)(b-a)^2+2m^2(m-1)(b-a)^2+m^2\mbox{max}\{a^2,b^2\}\nonumber\\
&{}&+4(1+\sqrt{2})^2m^2(m-1)^2(b-a)^2+2m^2(n-m)(b-a)^2\nonumber\\
&{}&+4(1+\sqrt{2})^2m^2(n-m)(m-1)(b-a)^2\nonumber\\
&:=&P(m,n)(b-a)^2+m^2\,\mbox{max}\{a^2,b^2\},
\end{eqnarray}
where
\begin{equation}\label{6.25}
P(m,n)=2m^2(m+n-2)+4(1+\sqrt{2})^2\,m^2(m-1)(n-1).
\end{equation}
\hfill$q.e.d$

\begin{Lemma}\label{Lemma6.6}
For the term \emph{(III)} , we have
\begin{eqnarray*}
\emph{(III)}\geq (m+2)\,b\,S-P(m)(b-a)S-4S^2,
\end{eqnarray*}
where $P(m)=2(m+2)+4(1+\sqrt{2})$.
\end{Lemma}

\emph{Proof}. Under an adapted unitary frame field, using the Gauss equation (\ref{7-6}), we have
\begin{eqnarray}\label{6.26}
\mbox{(III)}&=& \overline{a^\lambda_{i,j}}\big(a^\lambda_{p,j}\widetilde{R}^p_{ik\bar{k}}+a^\lambda_{i,p}\widetilde{R}^p_{jk\bar{k}}+2a^\lambda_{p,k}\widetilde{R}^p_{ij\bar{k}}-a^\mu_{i,j}\widetilde{R}^\lambda_{\mu k\bar{k}}\big)\nonumber\\
&{}& -\overline{a^\lambda_{i,j}}a^\lambda_{p,j}a^\mu_{i,k}\overline{a^\mu_{p,k}}
-\overline{a^\lambda_{i,j}}a^\lambda_{i,p}a^\mu_{j,k}\overline{a^\mu_{p,k}}
-2\overline{a^\lambda_{i,j}}a^\lambda_{p,k}a^\mu_{i,j}\overline{a^\mu_{p,k}}.
\end{eqnarray}
By using the bounds of holomorphic bisectional curvature, we have
\begin{eqnarray}\label{6.27}
\sum\limits_{i,j,k,p,\lambda}\overline{a^\lambda_{i,j}}a^\lambda_{p,j}\widetilde{R}^p_{ik\bar{k}}&=&\sum\limits_{j,k,\lambda}\sum\limits_{i,p}\overline{a^\lambda_{i,j}}a^\lambda_{p,j}\widetilde{R}^p_{ik\bar{k}}\nonumber\\
&\geq&a\sum\limits_{k}\sum\limits_{i,j,\lambda}|a^\lambda_{i,j}|^2=maS.
\end{eqnarray}
Similarly, we can get
\begin{equation}\label{6.28}
\sum\limits_{i,j,k,p,\lambda}\overline{a^\lambda_{i,j}}a^\lambda_{i,p}\widetilde{R}^p_{jk\bar{k}}\geq maS,\hspace{0.8cm}
\sum\limits_{i,j,k,\lambda,\mu}\overline{a^\lambda_{i,j}}a^\mu_{i,j}\widetilde{R}^\lambda_{\mu k\bar{k}}\leq mbS.
\end{equation}
On the other hand, by Lemma 6.4, the term
\begin{eqnarray}\label{6.29}
\sum\limits_{i,j,k,p,\lambda}\overline{a^\lambda_{i,j}}a^\lambda_{p,k}\widetilde{R}^p_{ij\bar{k}}&=&\sum\limits_{i,j,k,\lambda}\overline{a^\lambda_{i,j}}a^\lambda_{i,k}\widetilde{R}^i_{ij\bar{k}}+\sum\limits_{i,j,p,\lambda}\overline{a^\lambda_{i,j}}a^\lambda_{p,j}\widetilde{R}^p_{ij\bar{j}}\nonumber\\
&{}&-\sum\limits_{i,j,\lambda}\overline{a^\lambda_{i,j}}a^\lambda_{i,j}\widetilde{R}^i_{ij\bar{j}}+\sum\limits_{i\neq p,j\neq k,\lambda}\overline{a^\lambda_{i,j}}a^\lambda_{p,k}\widetilde{R}^p_{ij\bar{k}}\nonumber\\
&\geq& a \sum\limits_{i,\lambda}\sum\limits_{j}|a^\lambda_{i,j}|^2+a \sum\limits_{j,\lambda}\sum\limits_{i}|a^\lambda_{i,j}|^2\nonumber\\
&{}&-b\sum\limits_{i,j,\lambda}|a^\lambda_{i,j}|^2-2(1+\sqrt{2})(b-a)\sum\limits_{i\neq p,j\neq k,\lambda}|\overline{a^\lambda_{i,j}}|\,|a^\lambda_{p,k}|
\nonumber\\
&\geq& (2a-b)S-2(1+\sqrt{2})(b-a)S.
\end{eqnarray}

To estimate the remainning terms in (III), we set
\begin{equation*}
S'_{p\bar{i}}:=\sum\limits_{j,\lambda}a^\lambda_{p,j}\overline{a^\lambda_{i,j}},\hspace{0.6cm}
S''_{p\bar{j}}:=\sum\limits_{i,\lambda}a^\lambda_{i,p}\overline{a^\lambda_{i,j}}, \hspace{0.6cm}
S'''_{\lambda\,\bar{\mu}}:=\sum\limits_{i,j}a^\lambda_{i,j}\overline{a^\mu_{i,j}}.
\end{equation*}
It is clear that the matrices $S'=(S'_{p\bar{i}})$,  $S''=(S'_{p\bar{j}})$,  $S'''=(S'''_{\lambda\,\bar{\mu}})$ are Hermitian and semi-positive definite.
They are satisfy $S=\mbox{tr}(S')=\mbox{tr}(S'')=\mbox{tr}(S''')$. Locally, one can choose unitary frame to diagonalize $S'$, so we have
\begin{eqnarray}\label{6.30}
\sum\limits_{i,j,k,p,\lambda,\mu}\overline{a^\lambda_{i,j}}a^\lambda_{p,j}a^\mu_{i,k}\overline{a^\mu_{p,k}}
&=&\sum\limits_{i}S'_{i\bar{i}}S'_{i\bar{i}}\leq(\sum\limits_{i}S'_{i\bar{i}})^2=S^2.
\end{eqnarray}
Similarly, we also have
\begin{equation}\label{6.31}
\sum\limits_{i,j,k,p,\lambda,\mu}\overline{a^\lambda_{i,j}}a^\lambda_{p,k}a^\mu_{i,j}\overline{a^\mu_{p,k}}=\sum\limits_{\lambda}S'''_{\lambda\bar{\lambda}}S'''_{\lambda\bar{\lambda}}\leq(\sum\limits_{\lambda}S'''_{\lambda\bar{\lambda}})=S^2.
\end{equation}
By the trick of diagonalization, the term
\begin{eqnarray}\label{6.32}
\sum\limits_{i,j,k,p,\lambda,\mu}\overline{a^\lambda_{i,j}}a^\lambda_{i,p}a^\mu_{j,k}\overline{a^\mu_{p,k}}&=&\sum\limits_{j}S'_{j\bar{j}}S''_{j\bar{j}}\nonumber\\
&\leq& \big(\sum\limits_{j}(S'_{j\bar{j}})^2\big)^{1/2}\big(\sum\limits_{j}(S''_{j\bar{j}})^2\big)^{1/2}\leq S^2.
\end{eqnarray}
Thus, the lower bound of (III) follows from (\ref{6.26})-(\ref{6.32}).

\hfill$q.e.d$

Now we can prove the Simons integral inequality of pseudoholmorphic isometry. 

\begin{Theorem}
Let $(M,J,g)$ be a closed semi-K\"ahler manifold, and let $(\widetilde{M},\tilde{J},\tilde{g})$ be an almost Hermitian manifold pointwisely satisfies $a(x)\leq HB_{\widetilde{M}}\leq b(x)$ for $x\in \widetilde{M}$. Let $f:M\rightarrow \widetilde{M}$ be a pseudoholomorphic isometric immersion, then 
\begin{equation*}
\int_M\Big(\!-4S^2+\big[(m+2)b\!-\!P(m)(b-a)\big]S\!-\!P(m,n)(b-a)^2\!-m^2\,\emph{max}\{a^2,b^2\}\!\Big)dV_g\leq 0,
\end{equation*}
where $P(m)\!=\!2(m+2)+4(1+\sqrt{2})$, $P(m,n)\!=\!2m^2(m+n-2)+4(1+\sqrt{2})^2\,m^2(m-1)(n-1)$.
\end{Theorem}
 
\emph{Proof}. By Proposition 3.3, Lemma 6.5 and Lemma 6.6, we have
\begin{eqnarray}\label{6.33}
\frac{1}{2}\Box S&\geq& -\mbox{div}^c(X)-4S^2+\big[(m+2)b-P(m)(b-a)\big]S\nonumber\\
&{}&- P(m,n)(b-a)^2-m^2\,\mbox{max}\{a^2,b^2\}.
\end{eqnarray}
Notic that $(M,J,g)$ is closed and semi-K\"ahlerian, so $\int_M\,\Box S\,dV_g=0$ and $\int_M\,\mbox{div}^c(X)\,dV_g=0$, then the integral inequality follows from integrating on both sides of (\ref{6.33}).

\hfill$q.e.d$

\emph{Remark.} The semi-K\"ahlerian condition is to ensure that the terms $\Box S$ and $\mbox{div}^c(X)$ integrated on $M$ are zero. So, the same result holds when $(M,J,g)$ is quasi-K\"ahlerian, or almost-K\"ahlerian, or nearly-K\"ahlerian. 

\begin{Theorem}
Let $f:M\rightarrow \widetilde{M}$ be a pseudoholomorphic isometric immersion from almost Hermitian manifold $(M,J,g)$ into $(\widetilde{M},\tilde{J},\tilde{g})$ with parallel canonical second fundamental form. Suppose $(\widetilde{M},\tilde{J},\tilde{g})$ pointwisely satisfies $a(x)\leq HB_{\widetilde{M}}\leq b(x)$ for $x\in \widetilde{M}$, where $a(x)$ and $b(x)$ are bounded functions on $\widetilde{M}$. Then $S=0$, or 
$$\frac{1}{4}\big[(m+2)b_0-P(m)(b_0-a_0)\big]\leq S\leq \frac{m(n-m)}{m+n}\big[(m+2)a_0+P(m)(b_0-a_0)\big]$$ provided $P(m)(b_0-ma_0)< (m+2)b_0$, where  $P(m)=2(m+2)+4(1+\sqrt{2})$, $a_0:=\inf\limits_{x\in\widetilde{M}} a(x)$ and $b_{0}:=\sup\limits_{x\in \widetilde{M}} b(x)$.
\end{Theorem}

\emph{Proof}. Under the condition that $f$ has parallel canonical fundamental form, we have
\begin{equation}\label{6.34}
a^\alpha_{i,jk}=a^\alpha_{i,j\bar{k}}=0,
\end{equation}
and hence the Codazzi equation (\ref{7-11}) gives
\begin{equation}\label{6.35}
\widetilde{R}^\lambda_{ij\bar{k}}=0.
\end{equation}
The identity (\ref{4-4}) and (\ref{6.34}) imply that $S$ is a constant. For $X^k$ defined in (\ref{6.22}), the identity (\ref{6.35}) gives
\begin{equation}\label{6.36}
X^k=\overline{a^\lambda_{i,j}}\widetilde{R}^\lambda_{ij\bar{k}}=0.
\end{equation}
Thus, the term $\mbox{(II)}=0$ follows from (\ref{6.21}), (\ref{6.35}) and (\ref{6.36}). By using that $S$ is a constant, from Lemma 6.6, we have
\begin{equation}\label{6.37}
0\geq \big[(m+2)b_0-P(m)(b_0-a_0)-4S\big]S,
\end{equation} 
where $a_0:=\inf\limits_{x\in\widetilde{M}} a(x)$ and $b_{0}:=\sup\limits_{x\in \widetilde{M}} b(x)$.

On the other hand, modify the proof of Lemma 6.6, we can get the following estimates:
\begin{equation}\label{6.38}
\sum\limits_{i,j,k,p,\lambda}\overline{a^\lambda_{i,j}}a^\lambda_{p,j}\widetilde{R}^p_{ik\bar{k}}\leq mb_0S,
\end{equation}
\begin{equation}\label{6.39}
\sum\limits_{i,j,k,p,\lambda}\overline{a^\lambda_{i,j}}a^\lambda_{i,p}\widetilde{R}^p_{jk\bar{k}}\leq mb_0S,
\end{equation}
\begin{equation}\label{6.40}
\sum\limits_{i,j,k,\lambda,\mu}\overline{a^\lambda_{i,j}}a^\mu_{i,j}\widetilde{R}^\lambda_{\mu k\bar{k}}\geq ma_0S,
\end{equation}
\begin{eqnarray}\label{6.41}
\sum\limits_{i,j,k,p,\lambda}\overline{a^\lambda_{i,j}}a^\lambda_{p,k}\widetilde{R}^p_{ij\bar{k}}
\leq (2b_0-a_0)S+2(1+\sqrt{2})(b_0-a_0)S,
\end{eqnarray}
\begin{eqnarray}\label{6.42}
\sum\limits_{i,j,k,p,\lambda,\mu}\overline{a^\lambda_{i,j}}a^\lambda_{p,j}a^\mu_{i,k}\overline{a^\mu_{p,k}}
\geq\frac{1}{m}S^2,
\end{eqnarray}
\begin{equation}\label{6.43}
\sum\limits_{i,j,k,p,\lambda,\mu}\overline{a^\lambda_{i,j}}a^\lambda_{p,k}a^\mu_{i,j}\overline{a^\mu_{p,k}}\geq \frac{1}{n-m}S^2,
\end{equation}
\begin{eqnarray}\label{6.44}
\sum\limits_{i,j,k,p,\lambda,\mu}\overline{a^\lambda_{i,j}}a^\lambda_{i,p}a^\mu_{j,k}\overline{a^\mu_{p,k}}\geq 0.
\end{eqnarray}
From (\ref{6.38})-(\ref{6.44}) and Proposition 3.3, together with $S$ is a constant, we have
\begin{equation}\label{6.45}
0\leq \big[(m+2)a_0+P(m)(b_0-a_0)-\frac{n+m}{m(n-m)}S\big]S.
\end{equation}
Then, the bounds of $S$ follows from identities (\ref{6.37}) and (\ref{6.45}).

\hfill$q.e.d$

\hspace{0.6cm}

\noindent \textbf{Acknowledgments}. The second author would like to express his gratitude to professor X. Zhang for share the reference \cite{[Zh-12]}, and to X.Ch.Zhou for suggestions and comments in pareparing this paper. This project is supported by the NSFC (Grant No.11871445) and the Fundamental Research Funds for the Central Universities.


\end{document}